\documentclass[pdflatex,sn-mathphys-num]{sn-jnl}% Math and Physical Sciences Numbered Reference Style 
\usepackage{lmodern}
\usepackage{mathrsfs}%
\usepackage{graphicx}%
\usepackage{multirow}%
\usepackage{amsmath,amssymb,amsfonts}%
\usepackage{amsthm}%
\usepackage[title]{appendix}%
\usepackage{xcolor}%
\usepackage{textcomp}%
\usepackage{manyfoot}%
\usepackage{booktabs}%
\usepackage{algorithm}%
\usepackage{algorithmicx}%
\usepackage{algpseudocode}%
\usepackage{listings}%

\usepackage{pgfplots}
\pgfplotsset{compat=1.15}
\usepackage{mathrsfs}
\usetikzlibrary{arrows}
\definecolor{ccqqqq}{rgb}{0.8,0.,0.}
\definecolor{rvwvcq}{rgb}{0.08235294117647059,0.396078431372549,0.7529411764705882}
\definecolor{sexdts}{rgb}{0.1803921568627451,0.49019607843137253,0.19607843137254902}

\newtheorem{theorem}{Theorem}%  meant for continuous numbers
\newtheorem{proposition}[theorem]{Proposition}%
\newtheorem{remark}[theorem]{Remark}%

\newtheorem{definition}[theorem]{Definition}%
\newtheorem{lemma}[theorem]{Lemma}
\newtheorem{corollary}[theorem]{Corollary}

\newcommand{\cF}{\mathcal{F}}
\newcommand{\cI}{\mathcal{I}}
\newcommand{\R}{\mathbb{R}}
\newcommand{\N}{\mathbb{N}}
\newcommand{\cV}{\mathcal{V}}
\newcommand{\cO}{\mathcal{O}}
\newcommand{\cL}{\mathcal{L}}

\DeclareMathOperator*{\argmax}{arg\,max}
\DeclareMathOperator*{\argmin}{arg\,min}
\DeclareMathOperator{\cl}{cl}
\DeclareMathOperator{\conv}{conv}
\DeclareMathOperator{\diam}{diam}
\DeclareMathOperator{\dom}{dom}
\DeclareMathOperator{\relint}{relint}
\DeclareMathOperator{\range}{range}
\newcommand{\defeq}{\stackrel{\text{def}}{=}}

\begin{document}

\title[Article Title]{On the Acceleration of Proximal Bundle Methods}


\author[1]{\fnm{David} \sur{Fersztand}}\email{david282@mit.edu}

\author[1,2]{\fnm{Andy} \sur{Sun}}\email{sunx@mit.edu}

\affil[1]{\orgdiv{Operations Research Center}, \orgname{Massachusetts Institute of Technology}, \orgaddress{\street{77 Massachusetts Ave}, \city{Cambridge}, \postcode{02139}, \state{MA}, \country{USA}}}

\affil[2]{\orgdiv{Sloan School of Management}, \orgname{Massachusetts Institute of Technology}, \orgaddress{\street{77 Massachusetts Ave}, \city{Cambridge}, \postcode{02139}, \state{MA}, \country{USA}}}

\abstract{The proximal bundle method (PBM) is a fundamental and computationally effective algorithm for solving nonsmooth optimization problems. In this paper, we present the first variant of the PBM for smooth objectives, achieving an accelerated convergence rate of $\cO(\frac{1}{\sqrt{\epsilon}}\log(\frac{1}{\epsilon}))$, where $\epsilon$ is the desired accuracy. Our approach addresses an open question regarding the convergence guarantee of proximal bundle type methods, which was previously posed in two recent papers. We interpret the PBM as a proximal point algorithm and base our proposed algorithm on an accelerated inexact proximal point scheme. Our variant introduces a novel \textit{null step test} and oracle while maintaining the core structure of the original algorithm. The newly proposed oracle substitutes the traditional cutting planes with a smooth lower approximation of the true function. We show that this smooth interpolating lower model can be computed as a convex quadratic program. We also examine a second setting where Nesterov acceleration can be effectively applied, specifically when the objective is the sum of a smooth function and a piecewise linear one. }




\maketitle
% \tableofcontents

\section{Introduction}
\subsection{Problem setup}
We consider the Proximal Bundle type algorithms applied to the unconstrained minimization problem 
\begin{align}
    \underset{x\in \R^n}{\min}\, f(x), \label{pb: Initial pb}
\end{align}
where $f:\R^n\rightarrow \R$ is convex and $L_f$-smooth, i.e. $f$ has Lipschitz continuous gradient:
\begin{align*}
\|\nabla f(y) - \nabla f(x)\| \leq L_f \|x - y\|, \quad \forall x,y \in \R^n.
\end{align*}

\subsection{Motivation}
In \cite{LiangMonteiro2023}, Liang and Monteiro derived a convergence rate of $\cO(\frac{1}{\epsilon^2})$ up to logarithmic terms for the problem \eqref{pb: Initial pb}. They conclude their analysis by noting that this convergence rate \textit{is not optimal when $L_f > 0$. It would be interesting to design an accelerated variant of the PBM that is optimal when $L_f>0$}. 

In the same line of work, Diaz and Grimmer in \cite{Diaz2023optimal} studied the convergence rates of the classic PBM for a variety of non-smooth convex optimization problems. They concluded the paper with the statement that \textit{``it is likely that a variant of the bundle method can achieve accelerated rates, either by enforcing additional assumptions about the models $f_k$ or by
modifying the logic of the algorithm. We leave this as an intriguing open question for future research."}

In light of these remarks, it seems that the acceleration of the PBM for a smooth objective remains an open question. The complexity analysis of PBM is notoriously difficult. It is our hope that the investigation of modified bundle methods, such as the one we propose in Algorithm \ref{alg: accelerated prox bundle}, may ultimately lead to a better understanding of PBM. For smooth convex optimization problems, when the value of the smoothness parameter is available, there exist methods with a convergence rate $\cO(\frac{1}{\sqrt{\epsilon}})$, such as the accelerated gradient descent algorithm. This rate is known to be optimal for first-order methods \cite{drori_exact_2017}. We prove that, with a slight modification of the bundle minimization oracle and the \emph{null step test}, but keeping the overall structure of the algorithm, an accelerated rate can be achieved up to log terms. More precisely, an $\epsilon$-optimal solution can be reached in $\cO\left(\frac{1}{\sqrt{\epsilon}}\log\left(\frac{1}{\epsilon}\right)\right)$ iterations improving upon the standard PBM for which the best-known rate is $\cO(\frac{1}{\epsilon})$. 

When the objective is no longer smooth, but is the sum of a smooth convex function and a piecewise linear convex function, the best-known convergence rate is $\mathcal{O}\left(\frac{1}{\epsilon^{4/5}}\log\left( \frac{1}{\epsilon}\right)\right)$ \cite{2024proxbundlepolyhedral}. This paper proposes an accelerated variant of PBM for this setup that exhibits the same accelerated convergence rate of $\cO\left(\frac{1}{\sqrt{\epsilon}}\log\left(\frac{1}{\epsilon}\right)\right)$.

\subsection{Contribution}
\begin{enumerate}
    \item We establish the existence of a minimal smooth lower interpolating model, demonstrating that it can be constructed by solving a convex quadratically constrained quadratic program (QCQP). 
    \item We introduce a novel Accelerated Proximal Bundle algorithm, which incorporates modifications to both the oracle and the \textit{null step test} while preserving the fundamental structure and logic of the traditional proximal bundle algorithm. Moreover, the proposed \textit{null step test} does not rely on the tolerance $\epsilon$, positively answering an open question in \cite{LiangMonteiro2023}.
    \item We show that our proposed algorithm achieves an iteration complexity of $\cO\left(\frac{1}{\sqrt{\epsilon}}\log\left(\frac{1}{\epsilon}\right)\right)$. This result addresses an open question posed by \cite{LiangMonteiro2023} and \cite{Diaz2023optimal}, providing the first accelerated rate for the PBM.
    \item We propose an acceleration scheme for composite objectives, consisting of a smooth component and a piecewise linear one, and establish an iteration complexity of $\mathcal{O}\left(\frac{1}{\sqrt{\epsilon}} \log\left(\frac{1}{\epsilon}\right)\right)$.
\end{enumerate}

\subsection{Related work}
\paragraph{Kelley's method and bundle methods}
Kelley’s cutting plane method (1960) \cite{Kelley1960} is a foundational algorithm for minimizing non-smooth convex objectives. It iteratively refines a piecewise linear lower approximation of the true function. The PBM \cite{Kiwiel1990,deOliveira2014}, the trust region bundle method \cite{Schramm1992}, and the level bundle method \cite{Lemarechal1995, Lan2015} are variants of Kelley's method. Bundle methods were initially designed for non-smooth optimization because solving the proximal bundle subproblem is considered computationally expensive. This expense is justifiable for nonsmooth problems due to their inherently higher complexity. All these variants have been proven to converge to an optimal solution for any choice of parameter, contrary to gradient descent and its accelerated versions that require a stepsize chosen at most inversely proportional to the smoothness level. Likewise, subgradient methods necessitate diminishing stepsize sequences. These less complex algorithms might not converge if the stepsizes are not cautiously handled, thereby offering a persuasive argument to contemplate bundle methods. 

\paragraph{Convergence rates of the proximal bundle algorithm}
In 2000, Kiwiel \cite{Kiwiel2000} established the first convergence rate for the proximal bundle method, proving that an $\epsilon$-minimizer can be found in $\cO(\frac{1}{\epsilon^3})$ iterations.  \cite{Diaz2023optimal} improved the complexity analysis in the unconstrained case while restricting the proximity parameter schedule. They provided convergence rates for all combinations of smoothness and strong convexity. Namely, $\mathcal{O}(\frac{1}{\epsilon^2})$ when the function is only Lipschtiz continuous, $\mathcal{O}(\frac{1}{\epsilon})$ when the function is smooth or strongly convex, removing the log term in \cite{Du2017}. When the function is smooth and strongly convex, the convergence rate becomes $\mathcal{O}(\log\left(\frac{1}{\epsilon}\right))$. In \cite{liang2021proximal} and \cite{LiangMonteiro2023}, the authors proposed a new type of \textit{null step test} that is discussed in Section \ref{section: main ideas}. Their analysis provides an $\mathcal{O}(\frac{1}{\epsilon^2}\log(\frac{1}{\epsilon}))$ convergence rate when $f$ is only Lipschitz continuous, $\mathcal{O}(\frac{1}{\epsilon}\log(\frac{1}{\epsilon}))$ when $f$ is smooth. 
\paragraph{Proximal point algorithm}
The PBM can also be interpreted as an inexact proximal point method \cite{Martinet1972, RockafellarProximal1976}. G\"{u}ler in \cite{Guler1992} first proposed an inexact accelerated version of the proximal point algorithm. It reaches an $\cO(\frac{1}{\sqrt{\epsilon}})$ convergence rate by requiring the approximation error of the proximal problem to be summable. Note that it does not require the objective function to be smooth. In \cite{SalzoVilla2012}, Salzo and Villa corrected a subtle error in the convergence proof in \cite{Guler1992} with a slight modification of the hypothesis. He and Yuan in \cite{He2012} proposed a practical inexact criterion for the proximal point algorithm, eligible for Nesterov-type acceleration.

\paragraph{Acceleration of bundle methods}
Lan \cite{Lan2015} showed that the bundle level method can be modified to incorporate Nesterov-type acceleration. The algorithm is agnostic about the smoothness of the function in the sense that it achieves the optimal convergence rate for both smooth and nonsmooth problems (respectively $\cO(\frac{1}{\sqrt{\epsilon}})$ and $\cO(\frac{1}{\epsilon^{2}})$) without requiring the smoothness of the parameter as an input. To the best of our knowledge, no similar acceleration scheme is available for the PBM.

\subsection{The Proximal Bundle and the Proximal Point Algorithms}
For the sake of completeness, we present the formulations of the classic PBM (Algorithm \ref{alg: prox bundle}) and an accelerated inexact proximal point algorithm (Algorithm \ref{alg: proximal point algorithm}). The two algorithms are presented for minimizing a proper, closed, nonsmooth convex objective $f$. 

\begin{algorithm}[H]
  \caption{Proximal Bundle Method (PBM)}
  \label{alg: prox bundle}
\begin{algorithmic}[1]
\Require $x_0 \in \R^n$, $\mathcal{I}_0$ an index set of initial cuts, $\rho>0$, $\beta \in (0,1)$
    \For{$k \geq 0$}
      \State Compute $y_{k+1}$ solving the following quadratic program 
      \begin{align}
         & \underset{t \in \mathbb{R}, y \in \mathbb{R}^n}{\min} \; t + \frac{\rho}{2} \|y - x_k\|^2 \label{pb:PBA bundle minimization step}\\
    \text{s.t.} \quad &\forall i \in \mathcal{I}_k \quad t - f(y_i) - \langle v_i, y - y_i\rangle 
    \geq 0 \nonumber   
      \end{align}
        \State Compute $f(y_{k+1})$ and $v_{k+1}\in\partial f(y_{k+1})$
      \State $\mathcal{I}_{k+1} \leftarrow \mathcal{I}_k \cup \{k+1\}$ \Comment{(update piecewise linear model)}
      \If {$\beta(f(x_k) - f_{k}(y_{k+1})) \leq f(x_k) - f(y_{k+1})$} \Comment{(null step test)}\label{Alg:PBA:nulsteptest}
        \State $x_{k+1} \gets y_{k+1}$ \Comment{\text{(serious step)}}
      \Else
        \State $x_{k+1} \gets x_k$\Comment{(null step)}
      \EndIf
    \EndFor
\end{algorithmic}
\end{algorithm}
At each iteration $k$, the classic PBM solves a Quadratic Program (QP) in \eqref{pb:PBA bundle minimization step} that minimizes the sum of a piecewise linear lower model $f_k(y):=\max_{i\in\mathcal{I}_k} f(y_i)+\langle v_i, y-y_i\rangle$ with subgradient $v_i\in\partial f(y_i)$ and a quadratic proximal term. The solution of this optimization is denoted by $y_{k+1}$. The proximal center, which will be denoted by the letter $x$ throughout the article, is only updated subject to the \emph{null step test} on line \ref{Alg:PBA:nulsteptest} of Algorithm \ref{alg: prox bundle}. Steps in which the proximal center is updated are called \emph{serious steps}, while those where it remains unchanged are referred to as \emph{null steps}. The proximal term and the \emph{null step test} stabilize the convergence dynamics of the PBM algorithm compared to the classic cutting plane algorithm.

We also give here an accelerated inexact proximal point algorithm (a-IPPA), which applies Nesterov acceleration to the usual PPA \cite{Guler1992,He2012}. 
For $x_j \in \R^n$ and a convex function $f$, the prox operator is defined as $\operatorname{prox}_{\rho, f}(x_j) = {\argmax_{y\in \R^n}}\, f(y) + \frac{\rho}{2}\|y - x_j\|^2$. 
\begin{algorithm}[H]
  \caption{Accelerated Inexact Proximal Point Algorithm (a-IPPA)}
  \label{alg: proximal point algorithm}
\begin{algorithmic}
\Require a convex $f$, $x_0 \in \R^n$, a sequence $\rho_j>0$
    \State Set $\zeta_0 = x_0$ and $t_0 = 1$
    \For{$j \geq 0$}
      \State $y_{j+1} \approx \operatorname{prox}_{\rho_j, f}(x_j)$ \Comment{(inexact proximal step)}
      \State $v_{j+1} \in \partial f(y_{j+1})$
    \State $t_{j+1} \gets \frac{1+\sqrt{1 + 4t_j^2}}{2}$
    \State $\zeta_{j+1} \gets x_j - \frac{1}{\rho_j}v_{j+1}$
    \State $x_{j+1} \gets \zeta_{j+1} + \frac{t_j - 1}{t_{j+1}}(\zeta_{j+1} - \zeta_j)$
    \EndFor
\end{algorithmic}
\end{algorithm}
Interpreting the PBM as a proximal point algorithm, we develop an accelerated version of the PBM by specifying the inexact proximal step $y_{j+1} \approx \operatorname{prox}_{\rho, f}(x_j)$. We will detail both the method for finding an approximate solution and the criteria denoted $\approx$ that must be satisfied in Algorithm \ref{alg: proximal point algorithm} before updating the variables $(t, \zeta, x)$.
\paragraph{Paper overview}

Section \ref{section: Alg for smooth objective} focuses on the problem introduced in \eqref{pb: Initial pb} with a smooth objective function. This study is organized as follows. In Section \ref{section: Algorithm description}, we introduce our proposed algorithm. Section \ref{section: main ideas} outlines the key ideas behind the algorithm design and convergence proof. The new oracle, which constructs a smooth lower model, is detailed in Section \ref{section : Building a smooth model}. In Section \ref{Translation of the inexact criterion to our setup}, we show how the proposed \emph{null step test} ensures an accelerated convergence rate for the number of \emph{serious steps}. In Section \ref{Bound on the number of null steps}, we provide a bound on the number of consecutive \emph{null steps}. The main result of this paper is the convergence rate of our proposed algorithm and is provided in Section \ref{section: Overall convergence rate}. 

Section \ref{section: Alg for poly objective} presents an accelerated scheme in another setup in which the objective is not supposed to be smooth anymore. Instead, it is the sum of a convex smooth function and a convex piecewise linear one. 
\subsection{Definitions and notations}
We define $\|\cdot\|$ as the Euclidean norm and $\langle \cdot , \cdot \rangle$ the corresponding inner product in the $n$-dimensional Euclidean space $\R^n$. 

The relative interior, closure, and convex hull of a set $S\subseteq \R^n$ are denoted as $\relint(S)$, $\cl(S)$, and $\conv(S)$, respectively.  

The conjugate \cite{rockafellar1970convex} of a proper, convex function $h$ is defined as 
\begin{align*}
    h^*(u) = \sup_{x \in \R^n} \left\{\langle x, u \rangle - h(x)\right\}.
\end{align*}

For $h_1, h_2$, two real-valued functions of $\R^n$, ``$\leq$'' denotes the following property  
\begin{align*}
    h_1\leq h_2 \Leftrightarrow h_1(x)\leq h_2(x), \quad\forall x\in \R^n.
\end{align*}

We suppose that the problem \eqref{pb: Initial pb} has at least one optimal solution, denoted as $x^*$. This work studies the worst-case iteration complexity for finding an \emph{$\epsilon$-optimal} solution. Namely, the number of first-order oracle queries to find $x \in \R^n$ such that $f(x) - f(x^*)\leq \epsilon$. To simplify the expression of the convergence rates, we assume that $\epsilon<1$.
\section{Accelerated PBM for smooth objectives}
\label{section: Alg for smooth objective}
\subsection{Algorithm description}
\label{section: Algorithm description}
In this section, we introduce the accelerated Proximal Bundle Method (a-PBM). It builds upon the accelerated inexact proximal point algorithm (Algorithm \ref{alg: proximal point algorithm}). 
% We define the constant $C = \frac{2L_f}{\rho}\left( \frac{\sqrt{2L_f}}{\sqrt{\rho}} + 1 \right)$ involved in {null step test} (see Theorem \ref{thm: test of aPBA is sufficient for acceleration} for details). 
Although Algorithm \ref{alg: accelerated prox bundle} and its corresponding analysis are presented with a fixed parameter $\rho$, the analysis can be adapted to allow for a decreasing proximity parameter between consecutive \emph{null step} sequences.

\begin{algorithm}[hbt]
  \caption{Accelerated Proximal Bundle Method (a-PBM)}
  \label{alg: accelerated prox bundle}
\begin{algorithmic}[1]
\Require $f$ an $L_f$ smooth convex function, $x_0 \in \R^n$, $\mathcal{I}_0$ an index set of initial cuts, $\rho>0$, and $C = \frac{2L_f}{\rho}( \sqrt{2L_f/\rho} + 1)$.
    \State Set $\zeta_0 = x_0$ and $t_0 = 1$
    \For{$k \geq 0$}
      \State Compute $y_{k+1}$ solving the following convex QCQP :
      \begin{subequations}      
      \begin{align}
             & \underset{t \in \mathbb{R}, y \in \mathbb{R}^n}{\min} \; t + \frac{\rho}{2} \|y - x_k\|^2 \label{eq: bundle minimization alg smooth bundle method:obj}\\
    \text{s.t.} \quad &\forall i \in \mathcal{I}_k \quad t - f(y_i) - \langle \nabla f(y_i), y - y_i\rangle 
    \geq \frac{\rho^2}{2L_f}\|y - x_k + \frac{1}{\rho}\nabla f(y_i)\|^2. \label{eq: bundle minimization alg smooth bundle method:constr}
      \end{align}
      \end{subequations}\label{eq: bundle minimization alg smooth bundle method}
        \State Compute $f(y_{k+1})$ and $\nabla f(y_{k+1})$
      \State $\mathcal{I}_{k+1} \leftarrow \mathcal{I}_k \cup \{k+1\}$\Comment{(update piecewise linear model)}
      % \State \textbf{//Null step test} : 
      \If{$C \|y_{k+1}-y_{k}\| \leq\|x_k-y_{k+1}\|$} \Comment{(null step test)}
        \State $t_{k+1} \gets \frac{1+\sqrt{1 + 4t_k^2}}{2}$
        \State $\zeta_{k+1} \gets x_k - \frac{1}{\rho}\nabla f(y_{k+1})$ \Comment{(serious step)}
        \State $x_{k+1} \gets \zeta_{k+1} + \frac{t_k - 1}{t_{k+1}}(\zeta_{k+1} - \zeta_k)$
      \Else
        \State $x_{k+1} \gets x_k$
        \State $\zeta_{k+1} \gets \zeta_k$ \Comment{(null step)}
        \State $t_{k+1} \gets t_k$
      \EndIf
      % \State $(\hat{v}_{k+1}, \hat{b}_{k+1}) \in \underset{(v, b) \in \mathcal{V}}{\argmax}\; (v^{T} y_{k+1} + b)$
    \EndFor
\end{algorithmic}
\end{algorithm}

\subsection{Main ideas behind our proposed algorithm} \label{section: main ideas}
\paragraph{An accelerated inexact proximal point algorithm} In \cite{LiangMonteiro2023}, a variant of the PBM is presented that can be interpreted as an inexact proximal point algorithm (Algorithm \ref{alg: proximal point algorithm}). In 
 \cite{He2012}, He and Yuan present a Nesterov acceleration of the classical proximal point algorithm with a convenient approximation error for the proximal point oracle. 

\paragraph{A novel null step test} 
The usual \textit{null step test} for the PBM compares the expected progress in the objective value to the actual progress made via
\begin{align*}
    \beta(f(x_k) - f_{k}(y_{k+1})) \leq f(x_k) - f(y_{k+1}), 
\end{align*}
with $\beta$, a parameter in (0,1).
This criterion does not seem to translate easily into an inexact criterion from the proximal point algorithm literature. In \cite{liang2021proximal} and \cite{LiangMonteiro2023}, the \emph{null step test} is changed to 
\begin{align*}
    f(z_{k+1}) - f_{k}(y_{k+1}) \leq \frac{\epsilon}{2} + \frac{\rho}{2}\|y_{k+1} - x_k\|^2, 
\end{align*}
where $z_{k+1}$ represents the best iterate for the proximal problem.
Using this criterion that includes a quadratic term, the number of {serious steps} for their modified version of PBM  is bounded by $\rho\|x_0 - x^*\|^2/\epsilon$. This test makes it difficult to relate to the accelerated inexact proximal point method literature as the error term $\frac{\rho}{2}\|y_{k+1} - x_k\|^2$ may be large. However, imposing a stricter test will likely increase the provable upper bound on the number of consecutive {null steps}. Indeed, the number of consecutive {null steps} is shown to be $\cO(\log(\frac{1}{\epsilon}))$ through the geometric convergence of the quantity $m_{k+1} := f_{k}(y_{k+1}) + \frac{\rho}{2}\|y_{k+1} - x_k\|^2$. Although this linear convergence still holds with the oracle proposed in Section \ref{section : Building a smooth model}, it is not sufficient to prove the upper bound on the number of consecutive {null steps} with a stricter criterion like a summable approximation error, which is the one initially proposed in the first accelerated proximal point algorithms \cite{Guler1992}. 

In order to overcome this difficulty, we prove the geometric decrease of another quantity: the gap between the best value of the current proximal problem encountered so far and the last value of the bundle subproblem. This quantity $\xi$ is introduced in Section \ref{subsec : bound number of null step}. As the quantity $\frac{\rho}{2}\|y_{k+1} - x_k\|^2$ still appears in the study of the number of consecutive {null steps}, a way to obtain the geometric decrease of $\xi$ is to enforce a \emph{null step test} that relates this quantity to the distance between the two last iterates. 

From this result, we still need to show that such a test is sufficiently stringent to ensure the accelerated convergence rate of the corresponding inexact proximal point algorithm. It is indeed the case provided that the lower models $f_k$ are $L_f$-smooth. 

\paragraph{Smooth lower model functions} The cutting plane algorithm iteratively builds a lower approximation of $f$ as the pointwise maximum of the supporting hyperplanes at the previously queried points. This model function is nonsmooth. It is also the lowest convex function that interpolates $f$ in the gradients and function values at the previously queried points \cite{drori_optimal_2016}. Instead, if the lower models are smooth, the inexact criterion of \cite{He2012} translates into a comparison between $\|y_{k+1} - y_k
\|$, the distance between the two last iterates, and $\|y_{k+1}-x_k\|$ the distance from the last iterate to the proximal center. Provided that this ratio is small enough, the number of serious steps is guaranteed to be $\cO(\frac{1}{\sqrt{\epsilon}})$.

\begin{remark}[Extending the algorithm to a composite setting]
\label{remark: composite setting}It is possible to consider a composite minimization setting throughout the proof as in \cite{LiangMonteiro2023}. The problem becomes 
\begin{align*}
    \min_{x\in \R^n}\, f(x) + g(x),
\end{align*}
with $g$ a $L_g$-smooth convex function. 

To maintain clarity, we have excluded the detailed modifications to the algorithm and the convergence rate proof, as they are straightforward. For instance, the bundle minimization oracle \eqref{eq: bundle minimization alg smooth bundle method} would be replaced by 
\begin{align*}
    &\underset{t \in \mathbb{R}, y \in \mathbb{R}^n, u \in \mathbb{R}^n}{\min} \; t + g(y) + \frac{\rho}{2} \|y - x_k\|^2 \\
    \text{s.t.} \quad &\forall i \in \mathcal{I}_j \quad t - f(y_i) - \langle \nabla f(y_i), y - y_i\rangle \geq \frac{1}{2L_f}\|u - \nabla f(y_i)\|^2.
\end{align*}
The constant $C$ of the \textit{null step test} would be replaced by $\frac{2L_f}{\rho}\left( \frac{\sqrt{2(L_f+L_g)}}{\sqrt{\rho}} + 1 \right)$. In this setup, the convergence rate  $\cO(\frac{1}{\sqrt{\epsilon}}\log(\frac{1}{\epsilon}))$ still holds.
\end{remark}

\subsection{Building a smooth model}
\label{section : Building a smooth model}

In this section, we show how the bundle minimization step \eqref{eq: bundle minimization alg smooth bundle method} given in Algorithm \ref{alg: accelerated prox bundle} leads to lower models that display suitable properties for proving the accelerated convergence rate of the algorithm. It is a convex Quadratically Constrained Quadratic Program (QCQP) which can be efficiently solved. 

Drori and Teboulle in \cite{drori_optimal_2016} first interpreted Kelley's cutting plane method as iteratively minimizing the lowest convex function that interpolates the true function at the previous iterates. The method selects the next iterate to minimize the lowest convex function value. Additionally, this model function is not smooth, meaning it cannot truly represent the function \( f \). However, when \( f \) is known to be smooth, we demonstrate that the lowest smooth model functions can be constructed by solving a nonconvex QCQP. We also show that this nonconvex QCQP is equivalent to a more tractable convex QCQP.

We consider a sequence of consecutive null steps $j$ following the serious step $k$. At iteration $j$, we consider the bundle set $[k,j] \subset \mathcal{I}_j$. 
Compared to the classic PBM, we modify the bundle minimization step in order to build a sequence of approximation functions $f_j:\R^n\rightarrow\R$ that exhibit the following properties that for all $j$,
\begin{enumerate}
    \item $f_j$ is convex and $L_f$-smooth.
    \item $f_j$ interpolates $f$ at the previous iterates, i.e.,  for all $i \in \mathcal{I}_j,\,f_j(y_i) = f(y_i)$ and $\nabla f_j(y_i) = \nabla f(y_i)$.
    % \item $f_j$ is convex.
    \item $f(y_{j+1})\geq f_j(y_{j+1}) \geq f_{j-1}(y_{j+1})$, where $y_{j+1} = \argmin_{y\in \R^n}\,f_j(y) + \frac{\rho}{2}\|y-x_{j}\|^2$ is the result of the bundle minimization step \eqref{eq: bundle minimization alg smooth bundle method} in Algorithm \ref{alg: accelerated prox bundle}.
\end{enumerate}
The first property is useful to show that the \textit{null step test} given in Algorithm \ref{alg: accelerated prox bundle} implies the inexact criterion of \cite{He2012}. The other properties are used when deriving the bound on the number of consecutive null steps.

\begin{definition}
    % With $y_{i+1}$, denoting the result of the bundle minimization step \eqref{eq: bundle minimization alg smooth bundle method} in 
    For each iteration $j$ of Algorithm \ref{alg: accelerated prox bundle}, $\mathcal{F}_{L_f}(\cI_j)$ is defined to be the set of convex, $L_f$-smooth functions $\hat{f}: \R^n\rightarrow\R$ such that $\hat{f}$ interpolates the function values and gradients of $f$ at $(y_i)_{i\in\mathcal{I}_j}$, i.e., $\forall i\in \cI_j$, $\hat{f}(y_i) = f(y_i)$ and $\nabla  \hat{f}(y_i) = \nabla f(y_i)$.
% \begin{enumerate}
%     \item ,
%        \item $\forall i\in \cI_j, \quad \nabla  \hat{f}(y_i) = \nabla f(y_i)$.
       % \item $\hat{f}$ is convex and $L_f$-smooth. 
% \end{enumerate}
\end{definition}

\begin{remark}[Comparison with the recent approach by Florea and Nesterov \cite{florea2024optimallowerboundsmooth}]
During the final editing of our paper, we became aware of a paper from Florea and Nesterov studying smooth lower models
\cite{florea2024optimallowerboundsmooth}. Their contribution focuses on a different algorithm, specifically the Gradient Method with memory. Although their approach diverges somewhat from ours, the conclusions drawn in Section 2.3 of their paper align with those presented in our Section \ref{section : Building a smooth model}. Notably, they demonstrate that a smooth lower bound can be computed by solving a convex QCQP. 

To achieve this, Florea and Nesterov consider a family of functions indexed by $i\in \mathcal{I}_j$, $\phi_i: \R^n \times \R^n \rightarrow \R$ and $p:\R^n \rightarrow \R$ such that 
\begin{align}
    \forall (y,u)\in \R^n \times \R^n,\, \phi_i(y,u) &\defeq f(y_i) + \langle \nabla f(y_i), y -y_i\rangle + \frac{1}{2L_f}\|u - \nabla f(y_i)\|^2 \nonumber\\
    p(y) &\defeq \min_{g \in \R^n} \max_{i\in \mathcal{I}_j}\phi_i(y, g). \label{def: Nesterov and Florea p function}
\end{align}

To establish that $p$ is smooth, convex, and interpolates $f$ at the points $(y_i)_{i \in \mathcal{I}_j}$, Florea and Nesterov reformulate the inner maximization in \eqref{def: Nesterov and Florea p function} as a maximization over the simplex, leveraging strong duality and Danskin's theorem. In contrast, our approach begins with a nonconvex QCQP known to possess the desired properties outlined at the beginning of this section. We then demonstrate that this problem can be reduced to a more tractable convex QCQP.
\end{remark}

\subsubsection{A nonconvex QCQP with the desired properties}
We now introduce a nonconvex but explicit formulation for the lower bound of smooth functions interpolating $f$ at the previous iterates. In the next subsection, we will show that this nonconvex problem is equivalent to a convex QCQP with a formulation close to the usual proximal bundle minimization subproblem. More precisely, when $L_f=+\infty$ (\textit{i.e.} the function is not smooth), the two bundle minimization problems \eqref{pb:PBA bundle minimization step} and \eqref{eq: bundle minimization alg smooth bundle method} match exactly.

At each \emph{null step}, we consider the following problem 
\begin{subequations}\label{pb: smooth minimization nonconvex}
\begin{align}
    &\underset{t \in \mathbb{R}, y \in \mathbb{R}^n, u \in \mathbb{R}^n}{\min} \; t + \frac{\rho}{2} \|y - x_k\|^2 \label{pb: smooth minimization 2 constraints}\\
    \text{s.t.} \quad &\forall i \in \mathcal{I}_j, \quad t - f(y_i) - \langle \nabla f(y_i), y - y_i \rangle \geq \frac{1}{2L_f}\|u - \nabla f(y_i)\|^2,  \label{pb: smooth minimization 2 constraints. Constraint 1}  \\
    &\forall i \in \mathcal{I}_j, \quad f(y_i) - t - \langle u, y_i - y\rangle \geq \frac{1}{2L_f}\|u - \nabla f(y_i)\|^2. \label{pb: smooth minimization 2 constraints. Constraint 2}%\tag{\theequation b} 
\end{align}
\end{subequations}
Note that \eqref{pb: smooth minimization nonconvex} is a nonconvex QCQP due to the bilinear term $\langle u, y_i-y\rangle$ in \eqref{pb: smooth minimization 2 constraints. Constraint 2}.
\begin{lemma}
\label{lemma: rewrite QCQP with function variable}
    The above nonconvex QCQP \eqref{pb: smooth minimization nonconvex} is equivalent to
    \begin{align}
    \underset{y \in \mathbb{R}^n, \hat{f} \in \mathcal{F}_{L_f}(\cI_j)}{\min} \; \hat{f}(y) + \frac{\rho}{2} \|y - x_k\|^2 =  \underset{y \in \mathbb{R}^n}{\min}\left\{\frac{\rho}{2} \|y - x_k\|^2 + \underset{\hat{f} \in \mathcal{F}_{L_f}(\cI_j)}{\min} \; \hat{f}(y)\right\}. \label{smooth minimization problem over functions}
\end{align}
In particular, this shows that \eqref{pb: smooth minimization nonconvex} is always feasible, as $f \in \mathcal{F}_{L_f}(\cI_j) \neq \varnothing$.
\end{lemma}
\begin{proof}
Problem \eqref{smooth minimization problem over functions} is feasible, because $f$ itself is in $\mathcal{F}_{L_f}(\mathcal{I}_j)$. Then, according to Theorem 4 and Corollary 1 in \cite{taylor2016smooth}, for any $\hat{f}\in\mathcal{F}_{L_f}(\mathcal{I}_j)$ and $y\in\R^n$, $(\hat{f}(y), y, \nabla \hat{f}(y))$ is feasible for \eqref{pb: smooth minimization nonconvex}. Thus, the optimal objective value of \eqref{smooth minimization problem over functions} is no less than that of \eqref{pb: smooth minimization nonconvex}. This also shows that \eqref{pb: smooth minimization nonconvex} is feasible. Conversely, for any feasible solution $(t,y,u)$ of \eqref{pb: smooth minimization nonconvex}, by the same result in \cite{taylor2016smooth}, there exists an $L_f$-smooth convex function $\hat{f}$ such that  $\hat{f}(y) = t$, $\nabla \hat{f}(y) = u$, and for all $i\in \mathcal{I}_j$, $\hat{f}(y_i) = f(y_i)$ and $\nabla \hat{f}(y_i) =\nabla f(y_i)$, i.e., $\hat{f}\in\cF_{L_f}(\cI_j)$. 
% \begin{itemize}
%     \item $\hat{f}(y) = t$,
%     \item $\nabla \hat{f}(y) = u$,
%     \item $\forall i\in \mathcal{I}_j\quad \hat{f}(x_i) = f(x_i)$ and $\nabla \hat{f}(x_i) =\nabla f(x_i)$. 
% \end{itemize}
Thus, the optimal objective value of \eqref{pb: smooth minimization nonconvex} is no less than that of \eqref{smooth minimization problem over functions}. This proves the equivalence between \eqref{pb: smooth minimization nonconvex} and \eqref{smooth minimization problem over functions}. 
\end{proof}

Following Rockafellar's Convex Analysis book (Theorem 5.6) \cite{rockafellar1970convex}, we define the convex hull of a collection of functions.
\begin{definition}[Convex hull of a collection of functions]
    The convex hull of an arbitrary collection of functions $\{h_i : i \in \cI\}$ on $\mathbb{R}^n$ is denoted by $\conv \{h_i : i \in \cI\}$. It is the convex hull of the pointwise infimum of the collection. 
    %it is the function $h$ obtained via the convex hull of the union of the epigraphs of the functions $h_i$. It is the greatest convex function $h$  on $\mathbb{R}^n$ such that $h(x) \leq h_i(x)$, $\forall x \in \mathbb{R}^n, \forall i \in \cI$. %(not necessarily proper)

\end{definition}
 
\begin{lemma}
\label{lemma: f_min exists}
 Let $f_j$ be the closure of the convex hull of all functions in $\mathcal{F}_{L_f}(\cI_j)$, i.e., $f_j = \cl(\conv\{\hat{f} : \hat{f}\in\cF_{L_f}(\cI_j)\})$. Then $f_j$ is an element of $\mathcal{F}_{L_f}(\cI_j)$. Thus, $f_j=\conv\{\hat{f} : \hat{f}\in\cF_{L_f}(\cI_j)\}$. 
 % In  particular, it uniquely defines $f_j \in \mathcal{F}_{L_f}(\cI_j)$ such that  
Moreover, $f_j$ is the lowest element in $\cF_{L_f}(\cI_j)$, i.e., 
\begin{align*}
    \forall \hat{f} \in  \mathcal{F}_{L_f}(\cI_j), \, \forall x\in \R^n,\, f_j(x)\le \hat{f}(x).
\end{align*}
\end{lemma}
\begin{proof}
Let $E : x\mapsto \inf_{\hat{f}\in \mathcal{F}_{L_f}(\cI_j)}\hat{f}(x)$. Note that the conjugate $E^*=(\cl(\conv(E)))^*$ \cite[Corollary 12.1.1]{rockafellar1970convex}. Thus, $E^{**}=\cl(\conv(E))=f_j$. In the following, we study the properties of $E$, $E^*$, and $E^{**}$. First, notice that $\forall i\in \cI_j, E(y_i) = f(y_i)<+\infty$. Moreover, $\forall\hat{f}\in\cF_{L_f}(\cI_j), i\in \cI_j, y\in\R^n$, $\hat{f}(y)\ge \hat{f}(y_i)+\langle \nabla \hat{f}(y_i), y-y_i\rangle=f(y_i)+\langle \nabla f(y_i), y-y_i\rangle$. Thus, $\forall i\in \cI_j, y\in \R^n, E(y) \geq f(y_i) + \langle \nabla f(y_i), y - y_i\rangle >-\infty$. This shows that $E$ is proper. The conjugate $E^*$ of $E$ is given as  
\begin{align*}
    \forall w\in \R^n, \quad E^*(w) &= \underset{x\in \R^n}{\sup}\, \langle w, x\rangle -  E(x)\\
    &= \underset{x\in \R^n}{\sup}\, \langle w, x\rangle -  \inf_{\hat{f}\in \mathcal{F}_{L_f}(\cI_j)}\hat{f}(x)\\
    &= \underset{x\in \R^n}{\sup}\, \langle w, x\rangle +  \sup_{\hat{f}\in \mathcal{F}_{L_f}(\cI_j)}(-\hat{f}(x))\\
    &= \sup_{\hat{f}\in \mathcal{F}_{L_f}(\cI_j)} \underset{x\in \R^n}{\sup}\, \langle w, x\rangle - \hat{f}(x)\\
    &= \sup_{\hat{f}\in \mathcal{F}_{L_f}(\cI_j)}\, \hat{f}^*(w).
\end{align*}
As $\hat{f}$ is proper, convex, and $L_f$-smooth, $\hat{f}^*$ is proper convex and $1/L_f$-strongly convex. The pointwise supremum of convex functions is convex. So $\sup_{\hat{f}\in \mathcal{F}_{L_f}(\cI_j)}\, \hat{f}^*(w) - \frac{1}{2L_f}\|w\|^2$ is convex. This shows that $E^*$ is $1/L_f$-strongly convex. Thus, $E^{**}$ is $L_f$-smooth. Also, $E^{**}$ interpolates $f$ in function value and gradient at $(y_i)_{i\in \cI_j}$. This shows that $E^{**} \in \mathcal{F}_{L_f}(\cI_j)$ and $E^{**}$ is a lower bound of $\mathcal{F}_{L_f}(\cI_j)$.

\end{proof}
The previous lemma shows that we can unambiguously define $f_j$ as the minimum of $\mathcal{F}_{L_f}(\cI_j)$. 
%  This, combined with the optimality of $y_j = \argmin_{y\in \R^n} \,\{f_{j-1}(y) + \rho/2\|y -x_k\|^2\}$, and the strong convexity of $y\mapsto f_{j-1}(y) + \rho/2\|y -x_k\|^2$ yields :
% \begin{align*}
%     f_{j-1}(y_{j}) + \frac{\rho}{2}\|y_{j} - x_k\|^2 + \frac{\rho}{2}\|y_{j} - y_{j+1}\|^2 &\leq f_{j-1}(y_{j+1}) + \frac{\rho}{2}\|y_{j+1} - x_k\|^2 \\
%     &\leq f_{j}(y_{j+1}) + \frac{\rho}{2}\|y_{j+1} - x_k\|^2, 
% \end{align*}
% which is a key inequality for deriving a bound on the number of consecutive null steps in \cite{LiangMonteiro2023}. 
The following theorem shows that problem \eqref{pb: smooth minimization nonconvex} iteratively builds lower models of $f$ with desirable properties for establishing the accelerated convergence rate of Algorithm \ref{alg: accelerated prox bundle}. 
\begin{theorem}
    Problem \eqref{pb: smooth minimization nonconvex} defines a lower model function $f_j$ that is $L_f$-smooth, convex, and satisfies 
    \begin{align*}
        f(y_{j+1})\geq f_j(y_{j+1}) \geq f_{j-1}(y_{j+1}),
    \end{align*}
    with $y_{j+1} = \argmin_{y\in \R^n}\,f_j(y) + \frac{\rho}{2}\|y-x_k\|^2$.
\end{theorem}
\begin{proof}
Following Lemma \ref{lemma: rewrite QCQP with function variable}, \eqref{pb: smooth minimization nonconvex} is equivalent to \eqref{smooth minimization problem over functions}. We know from Lemma \ref{lemma: f_min exists} that $f_j$ is well defined and is minimal among all functions of $\mathcal{F}_{L_f}(\cI_j)$ for all $y$, hence
\begin{align*}
    \underset{y \in \mathbb{R}^n}{\min}\,\frac{\rho}{2} \|y - x_k\|^2 + \underset{\hat{f} \in \mathcal{F}_{L_f}(\cI_j)}{\min} \; \hat{f}(y) = \underset{y \in \mathbb{R}^n}{\min}\,\frac{\rho}{2} \|y - x_k\|^2 + f_j(y).
\end{align*}
This shows that $(f_j(y_{j+1}), y_{j+1}, \nabla f_j(y_{j+1}))$ is optimal for \eqref{pb: smooth minimization nonconvex}. 
%where $y_{j+1}$ is the optimum of the above problem, which uniquely exists due to the strong convexity of the objective function.

We have constructed $f_{j}$ such that $f_{j} \in  \mathcal{F}_{L_f}(\cI_j)$ and $\forall \hat{f} \in  \mathcal{F}_{L_f}(\cI_j)$, $f_{j} \leq \hat{f}$. In particular, as we suppose, during a sequence of null steps, that $\mathcal{I}_{j-1} \subset \mathcal{I}_{j}$, we deduce $f_{j} \in  \mathcal{F}_{L_f}(\cI_{j-1})$ ensuring that $f_{j} \geq f_{j-1}$. In addition, as $f \in \mathcal{F}_{L_f}(\cI_j)$, $f\ge f_j$.
\end{proof}

The following proposition and corollary show that the gradients of $f_j$ are in the convex hull of the gradients of $f$ at the previous iterates. This leads to a bound on the gradients of $f_j$ which will be useful for proving that the second constraint \eqref{pb: smooth minimization 2 constraints. Constraint 2} in the nonconvex QCQP \eqref{pb: smooth minimization nonconvex} is not necessary. 

\noindent
\begin{minipage}{0.54\textwidth}
\vspace{5pt}

Intuitively, taken any $\hat{f}\in\cF_{L_f}(\cI_j)$,  we want to construct a transformation $K$ on $\hat{f}$ such that $K(\hat{f})$ ``clips'' the gradients of $\hat{f}\in\mathcal{F}_{L_f}(\cI_j)$ to $\conv(\{\nabla f(y_i)\}_{i\in\cI_j})$. Figure 1 illustrates this operation. Suppose the bundle $\cI_j$ contains two points $y_1$ and $y_2$. Their corresponding cutting planes are dotted black lines. A smooth interpolating function $\hat{f}\in\cF_{L_f}(\cI_j)$ is given in bold dashed blue. $K$ clips its gradients and transforms $\hat{f}$ into the function represented by a solid red line. If we apply $K$ on $K(\hat{f})$, we get $K(\hat{f})$, i.e., $K(K(\hat{f}))=K(\hat{f})$. Thus, $K$ is an idempotent endomorphism on $\cF_{L_f}(\cI_j)$.  
\vspace{5pt}
\end{minipage}
\hfill
\begin{minipage}{0.44\textwidth}
\begin{tikzpicture}[line cap=round,line join=round,>=triangle 45]
\begin{axis}[
x=0.67cm,y=0.36cm,
axis lines=middle,
ymajorgrids=false,
xmajorgrids=false,
xmin=-3.7,
xmax=4,
ymin=-1,
ymax=8,
xtick=\empty,   
ytick=\empty,
clip=false]
\clip(-4.659862761720444,-3.5543824262996195) rectangle (5.248242862545619,10.702097370952169);
% Plotting the main x^2 function
\addplot [samples=50,domain=-5.0:5.0,line width=2.pt,blue, dashed] {x^2};

% Highlighting the part between y1 and y2
\addplot [samples=50,domain=-1:1.7292263309934688,line width=2.pt,ccqqqq] {x^2};

% Plotting y1 and y2 as dotted lines
\addplot [line width=1.pt,dotted,domain=-4.659862761720444:5.248242862545619] {(-0.5-1.*x)/0.5};
\addplot [line width=1.pt,dotted,domain=-4.659862761720444:5.248242862545619] {(-1.4951118519005668--1.7292263309934688*x)/0.5};

% Additional lines for context
\addplot [line width=2.pt,color=ccqqqq,domain=-4.659862761720444:-1.0] {(--2.90476707716136--5.80953415432272*x)/-2.90476707716136};
\addplot [line width=2.pt,color=ccqqqq,domain=1.7292263309934688:5.248242862545619] {(-5.770370490111496--6.673933176574202*x)/1.9297454176342312};

% Drawing points for y1 and y2
\addplot [only marks, mark=*, fill=black] coordinates {(-1,1) (1.7292263309934688,2.9902237038011337)};
\draw[color=black] (-0.8,1.612777774188329) node {$y_1$};
\draw[color=black] (1.3,3.1751317245720863) node {$y_2$};

\draw[color=blue] (-2.3,8) node {\textbf{$\hat{f}$}};

\draw[color=ccqqqq] (-3.9,5) node {$K(\hat{f})$};

\draw[color=black] (3,-1.5) node 
  {Figure 1};
\end{axis}
\end{tikzpicture}
  \raggedleft
\end{minipage}%

The following proposition makes this precise.
\begin{proposition}
    We define a mapping $K$ on $\mathcal{F}_{L_f}(\cI_j)$ such that
    \begin{align*}
        K : \hat{f} &\mapsto \left(\hat{f}^* + \delta_{G_j}\right)^*, 
    \end{align*}
    where $\delta_{G_j}$ is the indicator function of $G_j:=\conv(\{\nabla f(y_i)\}_{i\in\cI_j})$, i.e.,
    \begin{align*}
    \delta_{G_j}(u) = 
    \begin{cases}
       0, & \text{if $u\in G_j$},\\
       +\infty,  & \text{otherwise}.
    \end{cases}
\end{align*}
    Then, $K(\hat{f})\in\cF_{L_f}(\cI_j)$ and 
    \begin{align*}
        \forall \hat{f}\in \mathcal{F}_{L_f}(\cI_j), \quad  K(\hat{f}) \leq \hat{f}.
    \end{align*}
\end{proposition}

\begin{proof}
Let $\hat{f}\in \mathcal{F}_{L_f}(\cI_j)$.
As $\hat{f}$ is $L_f$-smooth, $\hat{f}^*$ is $1/L_f$-strongly convex. So $K(\hat{f})^*$ is $1/L_f$-strongly convex, which implies that $K(\hat{f})$ is $L_f$-smooth. 
For $i \in \cI_j$, we have
\begin{align*}
    \partial K(\hat{f})^*(\nabla f(y_i)) &= \partial (\hat{f}^* + \delta_{G_j})(\nabla f(y_i))\\
    &= \partial \hat{f}^*(\nabla f(y_i)) + \partial \delta_{G_j}(\nabla f(y_i))\\ 
    &\supseteq \partial \hat{f}^*(\nabla f(y_i)) + \{0\}\\
    &= \partial \hat{f}^*(\nabla f(y_i)) \ni y_i.
\end{align*}
The second equality holds because $\relint(\dom(\hat{f}^*)) \cap \relint(G_j) \neq \varnothing$. The last containment is due to the fact that $x^*\in \partial f(x) \iff x\in \partial f^*(x^*)$ for any closed convex function $f$ (see Theorem 23.5 in \cite{rockafellar1970convex}). Applying this fact again, we get $y_i\in \partial K(\hat{f})^*(\nabla f(y_i)) \iff \nabla f(y_i)\in \partial K(\hat{f})(y_i)$.
% This shows that $\nabla K(\hat{f})(y_i) = \partial (K(\hat{f})^*)^{-1}(y_i) \ni \nabla f(y_i)$. 
Since $K(\hat{f})$ is smooth, $\nabla f(y_i)=\nabla K(\hat{f})(y_i)$ for all $i\in\cI_j$. That is, $K(\hat{f})$ interpolates the gradients of $f$ at all the previous iterates. Moreover, we show that $K(\hat{f})$ also interpolates the function values of $f$ at the previous iterates. Indeed, since $K(\hat{f})$ is closed and $\nabla f(y_i)=\nabla K(\hat{f})(y_i)$, applying Theorem 23.5 in \cite{rockafellar1970convex}, we have 
\begin{align*}
    K(\hat{f})(y_i) &= \langle \nabla f(y_i), y_i\rangle - K(\hat{f})^*(\nabla f(y_i)) \\
    &=\langle \nabla f(y_i), y_i\rangle - \hat{f}^*(\nabla f(y_i))\\
    &= f(y_i).
\end{align*}
The second equality is due to the definition of $K(\hat{f})$ and $\nabla f(y_i)\in G_j$. The last equality again follows from Theorem 23.5 in \cite{rockafellar1970convex}. 
Thus, we have shown that $K(\hat{f}) \in \mathcal{F}_{L_f}(\cI_j)$. 
Also, as $K(\hat{f})^* \geq \hat{f}^*$, $K(\hat{f}) \leq \hat{f}$.
\end{proof}

\begin{corollary}
\label{corollary: fmin grad bounded}
The function $f_j$ has a gradient that takes values in the convex hull of the gradient of $f$ at the previous iterates: $\range(\nabla f_j) \subseteq \conv(\nabla f(y_i)_{i\in \cI_j})$.
\end{corollary}
\begin{proof}
    $K(f_j) \in \mathcal{F}_{L_f}(\cI_j)$ and $K(f_j) \leq f_j$. By minimality of $f_j$, $K(f_j) = f_j$.  $K(f_j)^* = f_j^*$ further implies that $\range(\nabla f_j) = \dom(f_j^*) \subseteq \conv(\nabla f(y_i))$.
\end{proof}

\subsubsection{Dropping the second line of constraints}

In this section, we show that the nonconvex constraints \eqref{pb: smooth minimization 2 constraints. Constraint 2} of problem \eqref{pb: smooth minimization nonconvex} are unnecessary and that \eqref{pb: smooth minimization nonconvex} is equivalent to the same problem without the second line of constraints as
\begin{subequations}\label{pb: smooth minimization 1 constraint}
    \begin{align}
    &\underset{t \in \mathbb{R}, y \in \mathbb{R}^n, u \in \mathbb{R}^n}{\min} \; t + \frac{\rho}{2} \|y - x_k\|^2 \stepcounter{equation}\\%\tag{\theequation}\\
    \text{s.t.} \quad &\forall i \in \mathcal{I}_j, \quad t - f(y_i) - \langle \nabla f(y_i), y - y_i \rangle \geq \frac{1}{2L_f}\|u - \nabla f(y_i)\|^2.
\end{align}
\end{subequations}
\begin{lemma}
\label{lemma: KKT conditions to derive value of u}
    In problem \eqref{pb: smooth minimization 1 constraint}, any optimal solution $(t^*, y^*, u^*)$ satisfies $u^* = \rho(x_k - y^*)$.
\end{lemma}
\begin{proof}

We write the Lagrangian of this problem: 
\begin{align*}
    \mathcal{L}(t,y,u,\lambda) &= t + \frac{\rho}{2} \|y - x_k\|^2 \\
    &\quad+ \sum_{i\in \mathcal{I}_j}{\lambda_i\left(-t  + f(y_i) + \langle \nabla f(y_i), y - y_i \rangle + \frac{1}{2L_f}\|u - \nabla f(y_i)\|^2\right)}.
\end{align*}
By the KKT conditions, for any optimal $(t^*, y^*, u^*)$, there exist $\lambda_i$'s such that
\begin{align}
    \begin{cases}
        0=1 - \sum_{i\in \mathcal{I}_j}{\lambda_i},\\
        0=\sum_{i\in \mathcal{I}_j}{\lambda_i\frac{1}{L}(u^* - \nabla f(y_i))}, \\
        0=\rho(y^* - x_k) + \sum_{i\in \mathcal{I}_j}{\lambda_i\nabla f(y_i)}. 
    \end{cases}\label{cases: KKT conditions to get u}
\end{align}
We conclude that $u^* = \rho(x_k - y^*)$. 
\end{proof}

Substituting $u$ with its value at optimality given by Lemma \ref{lemma: KKT conditions to derive value of u} in \eqref{pb: smooth minimization 1 constraint} gives the bundle minimization step \eqref{eq: bundle minimization alg smooth bundle method} in Algorithm \ref{alg: accelerated prox bundle}. Compared to the problem \eqref{pb:PBA bundle minimization step} solved in the usual PBM, it shares the same objective function,  number of variables and number of constraints. However, the constraints are convex quadratic rather than linear. 
\newline

We will now prove the equivalence between \eqref{eq: bundle minimization alg smooth bundle method} and \eqref{pb: smooth minimization nonconvex}. We consider problem \eqref{pb: smooth minimization 1 constraint} before optimizing over $y$,
\begin{subequations}\label{eq: before y}
\begin{align}
    &\underset{t \in \mathbb{R}, u \in \mathbb{R}^n}{\min} \; t + \frac{\rho}{2} \|y - x_k\|^2\\
    \text{s.t.} \quad &\forall i \in \mathcal{I}_j \quad t \geq  f(y_i) + \langle \nabla f(y_i), y - y_i\rangle  + \frac{1}{2L_f}\|u - \nabla f(y_i)\|^2. \label{constraint: constraint in 1}
\end{align}
\end{subequations}
The quadratic term in the objective only depends on $y$. We define $\Delta(y)$ as the set of optimal multipliers and $S(y)$ as the set of optimal solutions $(t,u)$, which is a singleton as the problem is strongly convex in $u$. In problem \eqref{constraint: constraint in 1}, the optimal value of $t$ is given by
\begin{align}
    t = \underset{u\in \R^n}{\min}\, \underset{i\in \cI_j}{\max}\, {f(y_i) + \langle \nabla f(y_i), y - y_i\rangle + \frac{1}{2L_f}\|u - \nabla f(y_i)\|^2}. \label{def: definition of tilde f}
\end{align}
This defines a function $\tilde{f}(y)$. This also shows that $\tilde{f}(y) \geq l(y):= \max_{i \in \cI_j}\, f(y_i) + \langle \nabla f(y_i), y - y_i\rangle$, where $f(y_i) + \langle \nabla f(y_i), y - y_i\rangle$ is the usual cutting plane. 

\begin{lemma}
\label{lemma: f tilde is differentiable}
    $\tilde{f}$ is differentiable with a gradient being the optimal value of $u$ in \eqref{eq: before y}.
\end{lemma}
\begin{proof}
$\tilde{f}(y)$ is defined as the result of a parametric convex problem in $y$. The set of optimal solutions is nonempty. The directional regularity condition (Definition 4.8 of \cite{bonnans_perturbation_2000} and Theorem 4.9) is satisfied. Indeed, the derivative with respect to $t$ of the constraint \eqref{constraint: constraint in 1} is $-1$ and the ambient set for $t$ is $\R$. Finally, the optimal solution $(t^*(y), u^*(y))$ of \eqref{def: definition of tilde f} is continuous in $y$ as the problem is strongly convex in $u$. According to Theorem 4.24 in \cite{bonnans_perturbation_2000}, the Hadamard directional derivative of $\tilde{f}$ is 
\begin{align*}
    \tilde{f}'(y, d) &= \underset{(t,u) \in S(y)}{\inf}\underset{\lambda \in \Delta(y)}{\sup}\, D_y \cL(t,y,u,\lambda) d \\
    &=\underset{\lambda \in \Delta(y)}{\sup}\, D_y \cL(t^*(y),y, u^*(y),\lambda) d \\
    &=\underset{\lambda \in \Delta(y)}{\sup}\, D_y \left[ t {+} \sum_{i \in \mathcal{I}_j} \lambda_i \left( {-}t {+} f(y_i) {+} \langle \nabla f(y_i), y - y_i\rangle {+} \frac{1}{2L_f} \|u {-} \nabla f(y_i)\|^2 \right) \right] d\\
    &=\underset{\lambda \in \Delta(y)}{\sup}\, \left(\sum_{i\in \mathcal{I}_j}\lambda_i \nabla f(y_i)\right) d\\
    &= \langle u, d\rangle.
\end{align*}
The last line follows from \eqref{cases: KKT conditions to get u}. 
This shows that $\tilde{f}$ is differentiable and that its gradient is the optimal value of $u$ in \eqref{eq: before y}.
\end{proof}

Recall that by Lemma \ref{lemma: f_min exists}, $f_j(y)=\min_{\hat{f}\in\cF_{L_f}(\cI_j)} \hat{f}(y)$ is the lowest smooth interpolating function in $\cF_{L_f}(\cI_j)$. 
% We define the usual cutting plane 
% \begin{align*}
%     l(y) = \max_{i \in \cI_j}\, f(y_i) + \langle \nabla f(y_i), y - y_i\rangle. 
% \end{align*}
Consider the function $H_j := f_j - \tilde{f}$. We know that
\begin{itemize}
    \item $H_j$ is differentiable as Lemma \ref{lemma: f tilde is differentiable} shows that $\tilde{f}$ is differentiable and $f_j\in \mathcal{F}_{L_f}(\cI_j)$ is smooth.
    \item $H_j\geq 0$ as \eqref{pb: smooth minimization 1 constraint} removes a constraint compared to \eqref{pb: smooth minimization nonconvex}.
\end{itemize} 

The following technical result (Proposition \ref{proposition: can find a sequence converging to the sup}) provides an elementary real analysis result. The proof can be found in Appendix \ref{sec: proof of prop 12}.
\begin{proposition}
\label{proposition: can find a sequence converging to the sup}
    Let $h: \R^n \rightarrow \R$ be an upper-bounded differentiable function and $\bar{h} = \sup_{x \in \R^n}\, h(x)$. For all $\delta >0$, there exists $x$ such that  
    \begin{align*}
        \begin{cases}
            \bar{h} - h(x) \leq \delta,\\
            \|\nabla h(x)\| \leq \delta.
        \end{cases}
    \end{align*}
\end{proposition}
% \begin{proof}
%     Let $x_\delta \in \R^n$ such that $\bar{h} - h(x_\delta) \leq \delta$. 
%     We consider the problem
%     \begin{align}
%         \max_{x \in \R^n}\, h(x) - \frac{\delta}{2}\|x - x_\delta\| \label{pb: proof for vanishing gradient sequence}.
%     \end{align}
%     As $h$ is upper-bounded, \eqref{pb: proof for vanishing gradient sequence} can be restricted to a compact feasible set. By continuity of the objective, this problem admits an optimal solution $x_\delta^*$. In particular, as $h(x_\delta^*) \geq h(x_\delta)$, $\|\nabla h(x_\delta^*)\|\leq \delta$ would conclude the proof. We assume the opposite.  $h$ is differentiable at $x_\delta^*$ so there exists $t>0$ such that $z = x_\delta^* + t\nabla h(x_\delta^*)$ satisfies:
%     \begin{align*}
%         h(z) = h(x_\delta^* + t\nabla h(x_\delta^*)) &\geq h(x_\delta^*) + \frac{3}{4} \left(t\nabla h(x_\delta^*)^T\right) \nabla h(x_\delta^*) \\
%         &= h(x_\delta^*) + \frac{3}{4}\|x_\delta^* - z\|\cdot \|\nabla h(x_\delta^*)\| \\
%        \implies\, h(z) - \frac{\delta}{2}\|z - x_\delta\| &\geq h(x_\delta^*) + \delta \frac{3}{4}\|x_\delta^* - z\| - \frac{\delta}{2}\|z - x_\delta\| \\
%        &\geq h(x_\delta^*) + \delta \frac{3}{4}\|x_\delta^* - z\| - \frac{\delta}{2}(\|x_\delta^* - z\| + \|x_\delta^* - x_\delta\|) \\
%        &\geq h(x_\delta^*)  - \frac{\delta}{2}\|x_\delta^* - x_\delta\| + \frac{\delta}{4}\|x_\delta^* - z\|.
%     \end{align*}
%     This contradicts the maximality of $x_\delta^*$ and concludes the proof.
% \end{proof}

\begin{theorem}
    The nonconvex QCQP \eqref{pb: smooth minimization nonconvex} is equivalent to the convex QCQP \eqref{eq: bundle minimization alg smooth bundle method}. 
\end{theorem}
\begin{proof}
We use Corollary \ref{corollary: fmin grad bounded} to show that $H_j$ is bounded.
\begin{align*}
    \forall y\in \R^n,\,f_j(y) &= \underset{i \in \cI_j}{\max}\,f(y_i) +  \langle \nabla f(y_i), y - y_i\rangle + \frac{1}{2L_f}\|\nabla f_j({y}) - \nabla f(y_i)\|^2\\
    &\leq  l(y) + \frac{1}{2L_f}\diam(\{\nabla f(y_{i})\}_{i \in \cI_j})^2.
\end{align*}
It holds that $H_j = f_j - \tilde{f} \leq f_j - l\leq \frac{1}{2L_f}\diam(\{\nabla f(y_{i})\}_{i \in \cI_j})^2$. This shows that $H_j$ is bounded from above. Let $\bar{H}_j$ denote its supremum over $\R^n$.

Following Proposition \ref{proposition: can find a sequence converging to the sup}, we know that for all $\eta>0$, there exists $\tilde{y} \in \R^n$ such that $\bar{H}_j - H_j(\tilde{y})\leq \eta/2$ and with $w = \nabla H_j(\tilde{y})$, satisfying
\begin{align*}
    \|w\|\leq \frac{\eta L_f}{2\diam(\{\nabla f(y_{i})\}_{i \in \cI_j})}.
\end{align*}

Let $\tilde{u} = \nabla f_j(\tilde{y}) = \nabla \tilde{f}(\tilde{y}) + w$ and $i_0 \in \cI_j$ such that
\begin{align*}
    f_j(\tilde{y}) = f(y_{i_0}) +  \langle \nabla f(y_{i_0}),\tilde{y} - y_{i_0}\rangle + \frac{1}{2L_f}\|\tilde{u} - \nabla f(y_{i_0})\|^2.
\end{align*}

This translates for $\tilde{f}$ into
\begin{align*}
    \tilde{f}(\tilde{y}) &= \underset{i \in \cI_j}{\max}\,f(y_i) +  \langle \nabla f(y_i), \tilde{y} - y_i\rangle  + \frac{1}{2L_f}\|\nabla \tilde{f}(\tilde{y}) - \nabla f(y_i)\|^2\\
    &\geq f(y_{i_0}) +  \langle \nabla f(y_{i_0}), \tilde{y} - y_{i_0}\rangle  + \frac{1}{2L_f}\|\tilde{u} - w - \nabla f(y_{i_0})\|^2\\
    &= f(y_{i_0}) +  \langle \nabla f(y_{i_0}), \tilde{y} - y_{i_0}\rangle  + \frac{1}{2L_f}\left(\|\tilde{u} - \nabla f(y_{i_0})\|^2 + \|w\|^2 {-} 2\langle w, \tilde{u} {-} \nabla f(y_{i_0})\rangle \right)\\
    &\geq f_j(\tilde{y}) + \frac{1}{2L_f}\left( - 2\langle w, \tilde{u} - \nabla f(y_{i_0})\rangle \right)\\
    &\geq f_j(\tilde{y}) - \frac{1}{L_f}\|w\|\cdot \|\tilde{u} - \nabla f(y_{i_0})\|\\
    &\geq f_j(\tilde{y}) - \frac{\eta}{2}.
\end{align*}
Thus, $\bar{H}_j \leq H_j(\tilde{y}) + \eta/2 = f_j(\tilde{y}) -\tilde{f}(\tilde{y}) + \eta/2\leq \eta/2 + \eta/2 = \eta$. As this is true for all $\eta>0$, $\sup_{y\in\R^n}\,H_j(y) = 0$. But we also know that $H_j\geq 0$. This shows that $H_j=0$ and proves that problems \eqref{pb: smooth minimization 1 constraint} and \eqref{pb: smooth minimization nonconvex} are equivalent. Lemma \ref{lemma: KKT conditions to derive value of u} shows the equivalence between \eqref{pb: smooth minimization 1 constraint} and \eqref{eq: bundle minimization alg smooth bundle method}.

\end{proof}

\subsection{Translation of the inexactness criterion to our setup}
\label{Translation of the inexact criterion to our setup}
Now, we translate the inexactness criterion originally proposed for the accelerated inexact proximal point algorithm  (a-IPPA, Algorithm \ref{alg: proximal point algorithm}) in \cite{He2012} to the setting of the accelerated proximal bundle method (a-PBM, Algorithm \ref{alg: accelerated prox bundle}). The idea is that the consecutive null steps \eqref{eq: bundle minimization alg smooth bundle method} after a serious step $k$ in a-PBM can be viewed as inexactly solving the proximal point problem $\operatorname{prox}_{\rho, f}(x_k)$ in a-IPPA until the inexactness criterion is satisfied, then the next serious step $k+1$ starts. We propose the following inexactness criterion for the proximal step 
   \begin{align}
        \left\langle \nabla f(y_{k+1}) - \nabla f(x_k -\frac{1}{\rho}\nabla f(y_{k+1})),  y_{k+1} - (x_k -\frac{1}{\rho}\nabla f(y_{k+1})) \right\rangle \leq \frac{1}{2\rho}\|\nabla f(y_{k+1})\|^2. \label{ineq: assumption for acceleration}
    \end{align}
The following theorem bounds the number of steps in a-IPPA when the inexactness criterion \eqref{ineq: assumption for acceleration} is satisfied by the inexact proximal step. This also gives the number of serious steps in a-PBM, when the serious steps satisfy the inexactness criterion  \eqref{ineq: assumption for acceleration}.
\begin{theorem}[Theorem 5.1 \cite{He2012}]
\label{thm: accelerated prox point he 1/k2}
Suppose that we apply the accelerated proximal point algorithm (a-IPPA, Algorithm \ref{alg: proximal point algorithm}) to problem \eqref{pb: Initial pb} with the inexact oracle for solving the proximal step problem satisfying condition \eqref{ineq: assumption for acceleration}, then
\begin{align*}
    f(\zeta_k) - f(x^*) \leq \frac{2\rho\|x_0 - x^*\|^2}{(k+1)^2},\quad \forall k\geq 1.
\end{align*}
An $\epsilon$-solution can be obtained after at most the following number of iterations
\begin{align*}
    \left\lceil\frac{\sqrt{2\rho}\|x_0-x^*\|}{\sqrt{\epsilon}}\right\rceil.
\end{align*}
\end{theorem}

% \begin{remark}
%     There seems to be a typo in \cite{He2012} as the theorem only holds for $k\geq 1$ (in their notations it should be $k\geq 2$ instead of $k\geq 1$. Indeed the result is false for the first iterate, as there is no condition relating $\rho$ and $L_f$). Another typo is in the algorithms, where they have replaced a $\nabla f$ by a $f$. 
% \end{remark}

\begin{remark}
    It is worth noting that G\"{u}ler in \cite{Guler1992} has first introduced an accelerated inexact proximal point algorithm. However, Salzo and villa in \cite{SalzoVilla2012} have shown that the proof is flawed and the result does not hold with the approximation criterion chosen by G\"{u}ler. In G\"{u}ler's setup, the correct convergence rate is $\cO(\frac{1}{\epsilon})$. Even though \cite{He2012} builds upon Guler's algorithm, and was published before \cite{SalzoVilla2012}, the convergence proof technique is quite different and does not suffer from the same flaw.  
\end{remark}

In order for a serious step in a-PBM to satisfy the inexactness criterion, we need to design a \emph{null step test} that implies the inexactness criterion. This is done by the following theorem.

\begin{theorem}[A sufficient condition for satisfying the inexact criterion]
\label{thm: test of aPBA is sufficient for acceleration}
    Let $k$ be a serious step in the execution of Algorithm \ref{alg: accelerated prox bundle} and let $x_k$ be the corresponding proximal center. Let $j > k$ be an index in the sequence of null steps following the serious step $k$. If 
    \begin{align}
        C \|y_{j+1}-y_{j}\| \leq\|x_k-y_{j+1}\| \quad\quad\quad\text{with} \quad\quad\quad C=\frac{2L_f}{\rho}\left( \frac{\sqrt{2 L_f}}{\sqrt{\rho}} + 1 \right), \label{eq: null step test for a-PBM}
    \end{align}
    then the inexactness criterion \eqref{ineq: assumption for acceleration} is satisfied.
\end{theorem}

\begin{proof}

We bound the left-hand side of \eqref{ineq: assumption for acceleration}:
\begin{align}
    &\left\langle \nabla f(y_{j+1}) - \nabla f(x_k -\frac{1}{\rho}\nabla f(y_{j+1})), y_{j+1} - (x_k -\frac{1}{\rho}\nabla f(y_{j+1})) \right\rangle \nonumber\\
    &\leq  \|\nabla f(y_{j+1}) - \nabla f(x_k -\frac{1}{\rho}\nabla f(y_{j+1}))\|\cdot\| y_{j+1} - (x_k -\frac{1}{\rho}\nabla f(y_{j+1})) \|\label{ineq: he criterion 1}\\
    &\leq L_f\|y_{j+1} - (x_k -\frac{1}{\rho}\nabla f(y_{j+1}))\|\cdot\| y_{j+1} - (x_k -\frac{1}{\rho}\nabla f(y_{j+1})) \|\label{ineq: he criterion 2} \\
    &=\frac{L_f}{\rho^2}\|\rho(y_{j+1} - x_k) +\nabla f(y_{j+1})\|^2. \nonumber
\end{align}
Eq. \eqref{ineq: he criterion 1} uses the Cauchy-Schwartz and Eq. \eqref{ineq: he criterion 2} uses the $L_f$-smoothness of $f$.

Note that because of the optimality of $y_{j+1}$ for the minimization step \eqref{eq: bundle minimization alg smooth bundle method},
\begin{align*}
    \nabla f_j (y_{j+1}) =- \rho(y_{j+1} - x_k).
\end{align*}
As $f_j$ is $L_f$-smooth and interpolates $f$ at $y_j$, for all $x \in \R^n$, we have
\begin{align*}
    \quad \|\nabla f(y_j) - \nabla f_j(x)\| = \|\nabla f_j(y_j) - \nabla f_j(x)\| \leq L_f\|x - y_j\|.
\end{align*}
Thus,
\begin{align*}
    \|\rho(y_{j+1} - x_k) +\nabla f(y_{j+1})\| &= \|-\nabla f_j (y_{j+1}) +\nabla f(y_{j+1})\|\\
    &\leq \|-\nabla f_j (y_{j+1}) +\nabla f(y_{j})\| + \|\nabla f(y_{j+1})-\nabla f(y_{j})\|\\
    &\leq 2L_f\|y_{j+1}-y_{j}\|.
\end{align*}
We now consider the right-hand side of \eqref{ineq: assumption for acceleration}, 
\begin{align*}
    \frac{1}{\sqrt{2\rho}}\|\nabla f(y_{j+1})\|&\geq \frac{1}{\sqrt{2\rho}}(-\|\nabla f_j (y_{j+1}) - \nabla f(y_{j+1})\| + \rho\|x_k-y_{j+1}\|)\\
    &\geq \frac{1}{\sqrt{2\rho}}\left( -2L_f\|y_{j+1}-y_{j}\| + \rho\|x_k-y_{j+1}\| \right).
\end{align*}
The following condition and $\rho\|x_k-y_{j+1}\| -2L_f\|y_{j+1}-y_{j}\|\geq 0$ imply \eqref{ineq: assumption for acceleration}:
\begin{align}
    \left(\frac{2\sqrt{L_f}}{\rho} + \frac{\sqrt{2}L_f}{\sqrt{\rho}}  \right)\|y_{j+1} - y_j\| &\leq \frac{\sqrt{\rho}}{\sqrt{2}}\|x_k-y_{j+1}\|, \nonumber\\
    \frac{2L_f}{\rho}\left( \frac{\sqrt{2 L_f}}{\sqrt{\rho}} + 1 \right)\|y_{j+1} - y_j\| &\leq \|x_k-y_{j+1}\|.
    \label{test: final test for inexact criterion}
\end{align}

%     4L_f^2\left( \frac{\sqrt{L_f+L_g}}{\rho} +\frac{1}{\sqrt{2\rho}} \right)^2 \|y_{j+1}-y_{j}\|^2 \leq \frac{\rho}{2}\|x_k-y_{j+1}\|^2\\
%     \frac{4L_f^2}{\rho^2}\left( \frac{\sqrt{2(L_f+L_g)}}{\sqrt{\rho}} + 1 \right)^2 \|y_{j+1}-y_{j}\|^2 \leq\|x_k-y_{j+1}\|^2 \label{test: final test for inexact criterion}.
% \end{align}
Eq. \eqref{test: final test for inexact criterion} is more stringent than $\rho\|x_k-y_{j+1}\| -2L_f\|y_{j+1}-y_{j}\|\leq 0$, showing that the criterion \eqref{test: final test for inexact criterion} alone suffices to imply \eqref{ineq: assumption for acceleration}. 
\end{proof} 

\begin{remark}
    In \cite{LiangMonteiro2023} and \cite{liang2021proximal} the authors considered a modified version of the PBM that has the drawback that the \textit{null step test} (hence the algorithm logic) depends on the value of $\epsilon$. Algorithm \ref{alg: accelerated prox bundle} (as the classic PBM \cite{Diaz2023optimal}) does not require the knowledge of the desired precision $\epsilon$. This makes the algorithm more versatile and reliable in diverse settings where the optimal precision might not be predetermined. 
\end{remark}

\subsection{Bound on the number of null steps}
\label{Bound on the number of null steps}

In \cite{LiangMonteiro2023}, when the objective function is smooth, the number of consecutive null steps is shown to be $\cO(\log(\frac{1}{\epsilon}))$. However, our setup is different in two ways:
\begin{enumerate}
    \item The update rule for the proximal center is given by the Nesterov acceleration scheme, so the bound involving distance to proximal centers is not valid anymore. 
    \item Our \emph{null step test} \eqref{eq: null step test for a-PBM} may be harder to satisfy than the one in \cite{LiangMonteiro2023}.
\end{enumerate}
The first point will be discussed in the next subsection. In Subsection \ref{subsec : bound number of null step} we will recall the main inequalities from \cite{LiangMonteiro2023} and introduce a new quantity $\xi_j$ which converges towards zero with a geometric rate. Finding an upper bound on the first term of this sequence, right after a serious step, and a lower bound on this sequence that is a power of $\epsilon$ suffices to show that there are $\cO(\log(\frac{1}{\epsilon}))$ consecutive null steps. 
\subsubsection{A Technical result: bound on the distance between the proximal center and the optimal solutions}

We first give a technical lemma (Lemma \ref{lm: recurrence inequality}) providing a bound derived from a recurrence inequality. The proof can be found in Appendix \ref{sec: proof of Lemma 18}.
\begin{lemma}
\label{lm: recurrence inequality}
    If a sequence satisfies the following recurrence inequality
    \begin{align*}
        \forall k\geq 0,\quad r_{k+1} \leq (1+\frac{2}{k+2}) r_k  + \frac{2C'}{k+2}, 
    \end{align*}
    where $C'>0$, then it holds that
    \begin{align*}
      \forall k\ge 1,\quad r_k \leq \left(e^2\cdot r_0+\frac{C'\pi^2e^{({3+\frac{\pi^2}{3}})}}{3}\right) \cdot k^2.
    \end{align*}
\end{lemma}

In the execution of Algorithm \ref{alg: accelerated prox bundle}, denote $j_k$ the overall iteration count for the $k$-th serious step. Then, we define $x_{j_k}^* := {\arg\min_{y\in \R^n}}\, f(y) + \frac{\rho}{2}\|y - x_{j_k}\|^2$. The next lemma bounds two key quantities for the analysis of the number of consecutive null steps.
\begin{lemma}
\label{lemma: bound on the distance to x^*}
    There exists a constant $C''>0$, depending only on $\rho$, $x^*$, $x_1$, $f(x^*)$, and $f(x_1)$ such that for all $k>0$,
    \begin{align}
        \|x_{j_{k}} - x^*\|  &\leq \frac{3 C''}{\epsilon} \quad \quad \text{and}\quad \quad \|x^* - x_{j_k}^*\| \leq \frac{6 C''}{\epsilon}.\label{ineq: bound on x_k - x^*}
    \end{align}
\end{lemma}

\begin{proof}
Following Lemma 5.1 in \cite{He2012}, considering consecutive serious steps with index $j_k$ for $k=1,2,3,\dots$,
\begin{align*}
2t_{j_{k}}^2 v_{j_{k}} - 2t_{j_{k+1}}^2 v_{j_{k+1}} \geq \rho \|u_{j_{k+1}}\|^2 - {\rho} \|u_{j_{k}}\|^2 , \quad \forall k \geq 0,
\end{align*}
where \( v_{j_{k}} := f(\zeta_{j_{k+1}}) - f(x^*) \) and \( u_{j_{k}} := t_{j_{k}} \zeta_{j_{k+1}} - (t_{j_{k}} - 1)\zeta_{j_{k}} - x^* \).
By adding these inequalities, we have
\begin{align*}
    \rho \|u_{j_{k}}\|^2 &\leq  \rho \|u_0\|^2 + 2t_0^2 v_0 \leq \rho \|x_1 - x^*\|^2 + 2(f(x_1) - f(x^*)).
\end{align*}

We define $C' = \sqrt{\|x_1 - x^*\|^2 + \frac{2}{\rho}(f(x_1) - f(x^*))}$,  which is an upper bound on $\|u_{j_{k}}\|$ for all $k$. Using the definition of $u_{j_{k}}$, we get
\begin{align*}
     \|u_{j_{k}}\| = \|t_{j_{k}} (\zeta_{j_{k+1}} - \zeta_{j_{k}}) + ( \zeta_{j_{k}} - x^*)\| \geq  t_{j_{k}}\| (\zeta_{j_{k+1}} - \zeta_{j_{k}})\| - \| \zeta_{j_{k}} - x^*\|.
\end{align*}
This provides an upper bound on $\| \zeta_{j_{k+1}} - \zeta_{j_{k}}\|$ as
\begin{align*}
    \| \zeta_{j_{k+1}} - \zeta_{j_{k}}\| \leq \frac{1}{t_{j_{k}}}(\| \zeta_{j_{k}} - x^*\| + C').
\end{align*}

From this inequality, we get
\begin{align*}
    \|\zeta_{j_{k+1}} - x^*\|&\leq \|\zeta_{j_k} - x^*\|+\|\zeta_{j_{k+1}} - \zeta_{j_k}\| 
    \leq \|\zeta_{j_k} - x^*\|+\frac{1}{t_{j_{k}}}(\| \zeta_{j_{k}} - x^*\| + C')\\
    &= (1+\frac{1}{t_{j_{k}}}) \|\zeta_{j_k} - x^*\|  + \frac{C'}{t_{j_{k}}}.
\end{align*}
We also have the following classic result about the Nesterov acceleration scheme, namely, $t_{j_{k}} \geq \frac{k+2}{2}$, with $k$ counting serious steps only (see, for instance, Lemma 5.3 in \cite{He2012} or \cite{BeckFISTA}). This is a non-homogeneous linear recurrence for which we can derive a bound $\|\zeta_{j_{k+1}} - x^*\|\leq C''\cdot \frac{1}{\epsilon}$ as follows. 

Using Theorem \ref{thm: accelerated prox point he 1/k2}, an $\epsilon$-optimal solution is reached after at most $\left\lceil \frac{\sqrt{2\rho}\|x_0 - x^*\|}{\sqrt{\epsilon}}\right\rceil$ serious iterations. We apply the result of Lemma \ref{lm: recurrence inequality} with $r_0 = \|\zeta_0 - x^*\| = \|x_0 - x^*\|$ to conclude that at the $k^{th}$ serious step,
\begin{align*}
     \|\zeta_{j_{k}} - x^*\| 
    &\leq \left(e^2\|x_0 - x^*\|+\frac{C'\pi^2e^{({3+\frac{\pi^2}{3}})}}{3}\right) \cdot k^2\\
    &\leq \left(e^2\|x_0 - x^*\|+\frac{C'\pi^2e^{({3+\frac{\pi^2}{3}})}}{3}\right) \cdot \left(1+ \frac{\sqrt{2\rho}\|x_0 - x^*\|}{\sqrt{\epsilon}}\right)^2 \\
    &\leq  \left(e^2\|x_0 - x^*\|+\frac{C'\pi^2 e^{({3+\frac{\pi^2}{3}})}}{3}\right)\cdot \left(2+ \frac{4\rho\|x_0 - x^*\|^2}{\epsilon}\right)
    = \cO\left(\frac{1}{\epsilon}\right).
\end{align*}
Let $C''$ be such that $\forall k\in\N^*, \, \|\zeta_{j_{k}} - x^*\| \leq C''/\epsilon$. If we suppose that $\epsilon\leq 1$, we can take
\begin{align*}
C'' = \left(2+ 4\rho\|x_0 - x^*\|^2\right){\cdot}\left[e^2\|x_0 - x^*\|{+}(\rho \|x_1 {-} x^*\|^2 {+} 2(f(x_1) {-} f(x^*)))\frac{\pi^2 e^{({3+\frac{\pi^2}{3}})}}{3}\right].
\end{align*}

We will conclude the proof by giving a bound on $\|x_{j_{k+1}} - x^*\|$. The first line uses the update of $x$ given in Algorithm \ref{alg: accelerated prox bundle}. 
\begin{align*}
    \|x_{j_{k+1}} - x^*\| &=\left\|\zeta_{j_{k+1}} - x^* + \frac{t_{j_{k}}-1}{t_{j_{k+1}}}\left((\zeta_{j_{k+1}} - x^*) - (\zeta_{k} - x^*)\right)\right\|\\
    &\leq \|\zeta_{j_{k+1}} - x^*\| + \frac{t_{j_{k}}-1}{t_{j_{k+1}}}(\|\zeta_{j_{k+1}} - x^*\| + \|\zeta_{k} - x^*\|)\\
    &\leq  (1+2\frac{t_{j_{k}}-1}{t_{j_{k+1}}})\frac{C''}{\epsilon} 
    \leq (1+\frac{4t_{j_{k}}}{1+\sqrt{1+4t_{j_{k}}^2}})\frac{C''}{\epsilon} \\
    &\leq (1+\frac{4t_{j_{k}}}{2t_{j_{k}}\sqrt{1+1/(4t_{j_{k}}^2)}})\frac{C''}{\epsilon} \tag{\theequation} \stepcounter{equation} \label{ineq: ratio of t smaller than 1}\\
    &\le (1+\frac{2}{\sqrt{1}})\frac{C''}{\epsilon}
    \leq 3\frac{C''}{\epsilon}. 
\end{align*}

By optimality of $x^*$ and $x_{j_k}^*$ for their respective problems,
\begin{align*}
    f(x^*) + \frac{\rho}{2}\|x_{j_k}^* -x_{j_k}\|^2 \leq f(x_{j_k}^*) + \frac{\rho}{2}\|x_{j_k}^* - x_{j_k}\|^2 \leq f(x^*) + \frac{\rho}{2}\|x^* -x_{j_k}\|^2.
\end{align*}
Thus, $\|x_{j_k}^* -x_{j_k}\| \leq \|x^* -x_{j_k}\|$. This leads to
\begin{align*}
    \|x^* - x_{j_k}^*\|&\leq \|x^* - x_{j_k}\| + \|x_{j_k}^* - x_{j_k}\|
    \leq  2\|x^* - x_{j_k}\|
    \leq 6 \frac{C''}{\epsilon}.
\end{align*}

\end{proof}

\subsubsection{Main result on the number of consecutive null steps}
\label{subsec : bound number of null step}
We consider a serious step $k$ followed by null steps $j$ with $j>k$. 

We define
\begin{align*}
    m_j = f_{j-1}(y_j) + \frac{\rho}{2}\|y_j - x_k\|^2.
\end{align*}
This is exactly the value of the bundle problem \eqref{eq: bundle minimization alg smooth bundle method} solved at iteration $j$. Let
\begin{align}
    \begin{cases}
        z_j = \underset{z \in \{y_{k+1}, y_{k+2}, ..., y_j\}}{\argmin} \{f(z) + \frac{\rho}{2}\|z-x_k\|^2\}\\
        \xi_j = f(z_{j}) + \frac{\rho}{2}\|z_{j}-x_k\|^2 - m_{j},      
    \end{cases}\label{eq: def of z_j and xi_j}
\end{align}
where $\xi_j$ represents the discrepancy between the best value of the \textit{null step} problem $\min_y \{f(y) + \frac{\rho}{2}\|y-x_k\|^2\}$ observed thus far and the most recent value of this problem when the true function $f$ is substituted with its model $f_{j-1}$.

\begin{lemma}
\label{lemma: linear convergence of xi}
Let $\tau \in [0,1)$, a constant such that
\begin{align}
    \frac{\tau}{1-\tau} \geq \frac{L_f}{\rho} + C^2. \label{ineq: bound on tau}
\end{align}
 Then, the following linear convergence rate on the quantity $\xi_j$ holds 
\begin{align}
   \forall j> k, \quad \quad \xi_{j+1}  &\leq \tau  \xi_j. 
\end{align}
\end{lemma}

\begin{remark}
    Contrary to the case where $f$ is only Lipschitz continuous in \cite{LiangMonteiro2023}, $\tau$ does not depend on $\epsilon$ and enables to derive a linear convergence of $\xi_j$, hence an $\cO(\log(\frac{1}{\epsilon}))$ number of consecutive null steps. 
\end{remark}

\begin{proof}[Proof of Lemma \ref{lemma: linear convergence of xi}]
% The following inequalities are valid whatever the null step test is.

 We define
\begin{align*}
    l_f(u,v) = f(v) + \langle \nabla f(v) , u-v\rangle,
\end{align*}
which corresponds to the supporting hyperplane of $f$ at $v$. 

As the bundle set, at step $j+1$, contains the cut made at step $j$ and that we suppose that $\mathcal{I}_j \subset \mathcal{I}_{j+1}$,
\begin{align}
    m_{j+1} &\geq l_f(y_{j+1}, y_j),\label{ineq: first ineq on m}\\
    m_{j+1} &\geq f_{j-1}(y_{j+1}) + \frac{\rho}{2}\|y_{j+1} - x_k\|^2 \geq m_j + \frac{\rho}{2}\|y_j - y_{j+1}\|^2. \label{ineq: second ineq on m}
\end{align}
For \eqref{ineq: first ineq on m}, we used the convexity of $f_j$ and the fact that $f_j(y_j) = f(y_j)$ and $\nabla f_j(y_j) = \nabla f(y_j)$. 
The second inequality \eqref{ineq: second ineq on m} uses the $\rho$-strong convexity of the bundle subproblem, the optimality of $y_j$ for the problem $\min_{y\in \R^n} \,\{f_{j-1}(y) + \rho/2\|y -x_k\|^2\}$ and $f_{j} \geq f_{j-1}$.  

Putting these two inequalities together, we get
\begin{align*}
    m_{j+1} &\geq  (1-\tau)\left(l_f(y_{j+1}, y_j) \right) + \tau (m_j + \frac{\rho}{2}\|y_j - y_{j+1}\|^2)\\
    &=\tau m_j  + (1-\tau)\left(l_f(y_{j+1}, y_j) + \frac{\tau \rho}{2(1-\tau)}\|y_j - y_{j+1}\|^2\right).
\end{align*}
Using the smoothness of $f$,  $l_f(y_{j+1}, y_j)  \geq f(y_{j+1}) - \frac{L_f}{2}\|y_j - y_{j+1}\|^2$. 
\begin{align}
    m_{j+1} -\tau m_j &\geq (1-\tau)f(y_{j+1}) + (1-\tau)\left(\frac{\tau \rho}{2(1-\tau)} - \frac{L_f}{2}\right)\|y_j - y_{j+1}\|^2. \label{ineq: final inequality on mj}
\end{align}
Using the definitions of $z_j$ and  $\xi_j$,
\begin{align}
    \xi_{j+1} &= f(z_{j+1}) + \frac{\rho}{2}\|z_{j+1}-x_k\|^2 - m_{j+1} \label{ineq: xi linear convergence 1}\\
    &\leq \tau (f(z_{j+1}) + \frac{\rho}{2}\|z_{j+1}-x_k\|^2) + (1-\tau)(f(y_{j+1}) + \frac{\rho}{2}\|y_{j+1}-x_k\|^2) -  m_{j+1} \label{ineq: xi linear convergence 2}\\
    &\leq \tau (f(z_{j+1}) + \frac{\rho}{2}\|z_{j+1}-x_k\|^2) + (1-\tau)(f(y_{j+1}) + \frac{\rho}{2}\|y_{j+1}-x_k\|^2)\nonumber\\
    &\quad - \tau m_j -(1-\tau)f(y_{j+1}) - (1-\tau)\left(\frac{\tau \rho}{2(1-\tau)} - \frac{L_f}{2}\right)\|y_j - y_{j+1}\|^2\label{ineq: xi linear convergence 3}\\
    &\leq\tau (f(z_{j+1}) + \frac{\rho}{2}\|z_{j+1}-x_k\|^2 - m_j) + \frac{\rho(1-\tau)}{2}\|y_{j+1}-x_k\|^2 \nonumber\\
    &\quad - (1-\tau)\left(\frac{\tau \rho}{2(1-\tau)} - \frac{L_f}{2}\right)\|y_j - y_{j+1}\|^2\label{ineq: xi linear convergence 4}\\
    &\leq \tau  \xi_j + \frac{1-\tau}{2}\left[{\rho}C^2-\left(\frac{\tau \rho}{(1-\tau)} - {L_f}\right)\right]\|y_j - y_{j+1}\|^2\label{ineq: xi linear convergence 5}\\
    &\leq \tau  \xi_j + \frac{1-\tau}{2}\underbrace{\left[{\rho}C^2 + {L_f}-\frac{\tau \rho}{(1-\tau)}\right]}_{\leq 0}\|y_j - y_{j+1}\|^2\label{ineq: xi linear convergence 6}\\
    &\leq \tau  \xi_j. \label{ineq: xi linear convergence 7}
\end{align}
In \eqref{ineq: xi linear convergence 1}, we rewrite the definition of $\xi$. The inequality \eqref{ineq: xi linear convergence 2} involves splitting the terms to make the quantities $\tau$ and $(1-\tau)$ appear. We also use the optimality of $z_{j+1}$ that ensures that $f(y_{j+1}) + \frac{\rho}{2}\|y_{j+1}-x_k\|^2 \geq f(z_{j+1}) + \frac{\rho}{2}\|z_{j+1}-x_k\|^2$.
In \eqref{ineq: xi linear convergence 3}, we directly apply inequality \eqref{ineq: final inequality on mj}.
In \eqref{ineq: xi linear convergence 4}, the terms $f(y_{j+1})$ cancel.
Eq. \eqref{ineq: xi linear convergence 5} is a key step in this proof as it leverages the \textit{null step test}.
In \eqref{ineq: xi linear convergence 6}, we use the definition of $z_{j+1}$.
Finally, we conclude with \eqref{ineq: xi linear convergence 7}, where $\tau$ is chosen to be sufficiently close to 1 in accordance with \eqref{ineq: bound on tau}, ensuring that the resulting quantity is less than or equal to $0$.

This shows the linear convergence of the quantity $\xi_j$. 

\end{proof}

After proving the geometric decrease of $\xi$, the following lemma gives a lower bound on this quantity. Supposing no serious step occurs, we know the \textit{null step test} is not satisfied. If $\xi_j$ is small enough, then the iterates $y_{j}$ will be close to the proximal center $x_k$. The proximal problem becomes approximately the initial problem \eqref{pb: Initial pb} and it is possible to show that an $\epsilon$-optimal solution to \eqref{pb: Initial pb} has been found. The lower bound on $\xi$, which is polynomial in $\epsilon$ (specifically, $\epsilon^3$), is sufficient for our proof. This quantity will appear as an argument of the logarithm function in the convergence rate of the algorithm.

\begin{lemma}
\label{lemma: minimum on xi}
There exists a constant $A>0$ such that, if at iteration $j$, $\xi_j \leq A\epsilon^{3}$ and the \textit{null step} criterion is not satisfied, then $f(y_{j})- f(x^*) \leq \epsilon$.
\end{lemma}

\begin{proof}
Let $\delta >0$. We consider an iteration $j$ such that $\xi_j \leq \delta$. Suppose that the \textit{null step} criterion is not satisfied, i.e., $C^2\|y_{j+1}-y_{j}\|^2 \geq\|x_k-y_{j+1}\|^2$. Also, using  \eqref{ineq: second ineq on m},
\begin{align}
    m_{j+1} \geq m_j + \frac{\rho}{2}\|y_{j+1} - y_j\|^2.\label{ineq: null step test ineq on mj}
\end{align}

We combine these inequalities to get 
\begin{align*}
    % m_{j+1}&\geq m_j + \frac{\rho}{2}\|y_{j+1} - y_j\|^2\\
     m_{j+1} - m_j &\geq \frac{\rho}{2C^2}\|x_k-y_{j+1}\|^2.
\end{align*}
This provides an upper-bound on $\|x_k-y_{j+1}\|^2$,
\begin{align}
    \|x_k-y_{j+1}\|^2 &\leq \frac{2C^2}{\rho}\left(m_{j+1} - m_j \right)\label{ineq: distance iterate to center 1}\\
    &=  \frac{2C^2}{\rho}\left(f_j(y_{j+1}) + \frac{\rho}{2}\|y_{j+1}-x_k\|^2 - m_j \right)\label{ineq: distance iterate to center 2}\\
    &\leq  \frac{2C^2}{\rho}\left(f_j(z_{j}) + \frac{\rho}{2}\|z_{j}-x_k\|^2 - m_j \right)\label{ineq: distance iterate to center 3}\\
    &\leq  \frac{2C^2}{\rho}\left(f(z_{j}) + \frac{\rho}{2}\|z_{j}-x_k\|^2 - m_j \right) = \frac{2C^2}{\rho}\xi_j\label{ineq: distance iterate to center 4}\\
    &\leq  \frac{2C^2}{\rho}\delta.\label{ineq: distance iterate to center 5}
\end{align}
Eq. \eqref{ineq: distance iterate to center 2} is given by definition of $m_{j+1}$. Eq. 
\eqref{ineq: distance iterate to center 3} uses the optimality of $y_{j+1}$ for \eqref{eq: bundle minimization alg smooth bundle method}.
Eq. \eqref{ineq: distance iterate to center 4} uses $f_j\leq f$ and the definition of $\xi_j$.
Eq. \eqref{ineq: distance iterate to center 5} concludes using the upper bound on $\xi_j$ supposed in this lemma.

We use this to bound $\|y_j - x_k\|^2$ as
\begin{align}
    \|y_j - x_k\|^2&\leq 2\|x_k-y_{j+1}\|^2 + 2\|y_j-y_{j+1}\|^2\\
    &\leq 2\frac{2C^2}{\rho}\delta + 2\frac{2}{\rho}(m_{j+1}-m_j) \label{ineq: using the null step test again 1}\\
    &\leq \frac{4C^2}{\rho}\delta + 4\delta/\rho \label{ineq: using the null step test again 2}\\
    &=4\rho^{-1}\left(C^2 + 1 \right) \delta.\label{ineq: using the null step test again 3}
\end{align}
To derive \eqref{ineq: using the null step test again 1}, we use \eqref{ineq: null step test ineq on mj} again. For \eqref{ineq: using the null step test again 2}, $\forall z\in \R^n, \, m_{j+1}\leq f(z) + \frac{\rho}{2}\|z-x_k\|^2$. In particular, $m_{j+1} - m_j \leq f(z_{j}) + \frac{\rho}{2}\|z_{j}-x_k\|^2 - m_{j} = \xi_j\leq \delta.$

We want to relate the best value for the \textit{null step} problem seen so far to the current value of this problem, namely $f(z_j) + \frac{\rho}{2}\|z_j - x_k\|^2$ and $f(y_j) + \frac{\rho}{2}\|y_j - x_k\|^2$.
\begin{align}
    &f(z_j) + \frac{\rho}{2}\|z_j - x_k\|^2 \geq f_j(z_j) + \frac{\rho}{2}\|z_j - x_k\|^2 \label{ineq: replace yj by zj 1}\\
    &\geq  f_j(y_{j+1}) + \frac{\rho}{2}\|y_{j+1} - x_k\|^2 \label{ineq: replace yj by zj 2}\\
    &\geq f_j(y_j) + \langle \nabla f_j(y_j), y_{j+1} - y_j\rangle+ \frac{\rho}{2}\|y_{j+1} - x_k\|^2\label{ineq: replace yj by zj 3}\\
    &\geq f(y_j) - \| \nabla f(y_j)\| \cdot \| y_{j+1} - y_j\|\label{ineq: replace yj by zj 4}\\
    &=  f(y_j) - \| \nabla f(x^*) - \nabla f(y_j)\|\cdot \| y_{j+1} - y_j\|\label{ineq: replace yj by zj 5}\\
    &\geq  f(y_j) - L_f\|x^* -y_j\|\cdot \| y_{j+1} - y_j\|\label{ineq: replace yj by zj 6}\\
    &\geq  f(y_j) - L_f(\|x^* - x_k\|+\|y_j - x_k\|) \| y_{j+1} - y_j\|\label{ineq: replace yj by zj 7}\\
    &\geq  f(y_j) - L_f\left[3C''/\epsilon+4\rho^{-1}\left(C^2 + 1 \right) \delta\right]\cdot 2\delta/\rho + \left(\frac{\rho}{2}\|y_j - x_k\|^2 - 4\rho^{-1}\left(C^2 + 1 \right) \delta\right)\label{ineq: replace yj by zj 8}\\
    &=f(y_j) + \frac{\rho}{2}\|y_j - x_k\|^2 - 2\delta L_f/\rho\left[3C''/\epsilon+4\rho^{-1}\left(C^2 + 1 \right) \delta + 2\left(C^2 + 1 \right)/L_f\right]\label{ineq: replace yj by zj 9}\\
    &=f(y_j) + \frac{\rho}{2}\|y_j - x_k\|^2  - D\delta/\epsilon \label{ineq: replace yj by zj 10}.
\end{align}
Eq. \eqref{ineq: replace yj by zj 1} follows from the fact that $f_j \leq f$. Eq. \eqref{ineq: replace yj by zj 2} is derived from the optimality of $y_{j+1}$ in the bundle minimization problem described by \eqref{eq: bundle minimization alg smooth bundle method}. Eq. \eqref{ineq: replace yj by zj 3} utilizes the convexity of $f_j$. \eqref{ineq: replace yj by zj 4} applies the Cauchy-Schwarz inequality and the fact that $f_j$ and $f$ have the same value and gradient at $y_j$. Eq. \eqref{ineq: replace yj by zj 5} relies on the first-order optimality condition at $x^*$. Eq. \eqref{ineq: replace yj by zj 6} uses the smoothness of $f$. Eq. \eqref{ineq: replace yj by zj 7} employs the triangle inequality. Eq. \eqref{ineq: replace yj by zj 8} incorporates bounds from \eqref{ineq: using the null step test again 3}, \eqref{ineq: bound on x_k - x^*}, and \eqref{ineq: null step test ineq on mj}. Finally, \eqref{ineq: replace yj by zj 9} involves rearranging terms, and \eqref{ineq: replace yj by zj 10} introduces a new constant $D > 0$.

Using this inequality, we derive a bound on $\|x_k^* - y_j\|$ as
\begin{align}
    f(y_j) + \frac{\rho}{2}\|y_j - x_k\|^2 &\leq f(z_j) + \frac{\rho}{2}\|z_j - x_k\|^2 +D\delta/\epsilon \nonumber \\
    &= f(z_j) + \frac{\rho}{2}\|z_j - x_k\|^2  - m_j + m_j +D\delta/\epsilon \nonumber\\ 
    &= D\delta/\epsilon  + \xi_j + f_{j-1}(y_j) + \frac{\rho}{2}\|y_j - x_k\|^2 \label{ineq: xstark to yk 1}\\
    &\leq \delta+D\delta/\epsilon+ f_{j-1}(y_j) + \frac{\rho}{2}\|y_j - x_k\|^2\label{ineq: xstark to yk 2}\\
    &\leq \delta + D\delta/\epsilon + f_{j-1}(x_k^*) + \frac{\rho}{2}\|x_k^* - x_k\|^2 - \frac{\rho}{2}\|x_k^* - y_j\|^2\label{ineq: xstark to yk 3}\\
    &\leq \delta+D\delta/\epsilon+f(x_k^*) + \frac{\rho}{2}\|x_k^* - x_k\|^2 - \frac{\rho}{2}\|x_k^* - y_j\|^2\label{ineq: xstark to yk 4}\\
    &\leq \delta+D\delta/\epsilon+f(y_j) + \frac{\rho}{2}\|y_j - x_k\|^2 - {\rho}\|x_k^* - y_j\|^2\label{ineq: xstark to yk 5}\\
    \iff \|x_k^* - y_j\|^2 &\leq \delta(1+D/\epsilon)/\rho.\label{ineq: xstark to yk 6}
\end{align}
Eq. \eqref{ineq: xstark to yk 1} employs the definition of $\xi_j$, while \eqref{ineq: xstark to yk 2} utilizes the condition $\xi_j \leq \delta$. \eqref{ineq: xstark to yk 3} leverages both the optimality of $y_j$ in solving the bundle subproblem at iteration $j$ and the strong convexity of the corresponding function, whereas \eqref{ineq: xstark to yk 4} relies on the inequality $f_{j-1} \leq f$. Eq. \eqref{ineq: xstark to yk 5} depends on the optimality of $x_k^*$ for the proximal problem centered at $x_k$, alongside the strong convexity of the corresponding function. In \eqref{ineq: xstark to yk 6}, we rearrange the terms.

Starting from \eqref{ineq: xstark to yk 4} we get
\begin{align}
\delta(1+D/\epsilon)&\geq f(y_{j}) + \frac{\rho}{2}\|y_{j} - x_k\|^2 - (f(x_k^*) + \frac{\rho}{2}\|x_k^*- x_k\|^2 - \frac{\rho}{2}\|y_{j} - x_k^*\|^2) \nonumber\\
&\geq f(y_{j}) + \frac{\rho}{2}\|y_{j} - x_k\|^2 \nonumber\\
&\quad - \left(f(x^*) + \frac{\rho}{2}\|x^*- x_k\|^2 - \frac{\rho}{2}\|y_{j} - x_k^*\|^2 - \frac{\rho}{2}\|x^* - x_k^*\|^2\right).\label{ineq: towards global optimality 2}
\end{align}

Eq. \eqref{ineq: towards global optimality 2} uses the optimality of $x_k^*$ for the proximal subproblem and the strong convexity of this subproblem.

Rearanging, using $f_j\leq f$ and writing $x^* - x_k = (x^* - x_k^*) + (x_k^* - y_j) + (y_j - x_k)$,
\begin{align}
    &f(y_{j})- f(x^*)\nonumber \\&\leq \delta(1+D/\epsilon) - \frac{\rho}{2}\|y_{j} - x_k\|^2  + \frac{\rho}{2}\|x^*- x_k\|^2 - \frac{\rho}{2}\|y_{j} - x_k^*\|^2 - \frac{\rho}{2}\|x^* - x_k^*\|^2 \nonumber\\
    &\leq \delta(1{+}D/\epsilon) {+} \rho\left(\|y_{j} {-} x_k^*\|{\cdot}\|x^* {-} x_k^*\| {+} \|y_{j} {-} x_k\|{\cdot}\|x^* {-} x_k^*\| {+} \|y_{j} {-} x_k\|{\cdot}\|y_{j} {-} x_k^*\| \right) \label{ineq: small xi leads to optimality 1}\\
    &\leq \delta(1{+}D/\epsilon) {+} \delta(1{+}D/\epsilon)\|x^* {-} x_k^*\| {+} 4\left(C^2 {+} 1 \right) \delta\|x^* {-} x_k^*\| {+} 4\delta^2\left(C^2 {+} 1 \right) (1{+}D/\epsilon) \label{ineq: small xi leads to optimality 2} \\
       &= \delta\left[(1+D/\epsilon) + \left(4C^2 + 5 + D/\epsilon\right) \|x^* - x_k^*\| + 4\delta\left(C^2 + 1 \right) (1+D/\epsilon) \right]. \label{ineq: small xi leads to optimality}
\end{align}
We use the Cauchy-Schwarz inequality to derive \eqref{ineq: small xi leads to optimality 1}. For \eqref{ineq: small xi leads to optimality 2}, we apply \eqref{ineq: using the null step test again 3} and \eqref{ineq: xstark to yk 6}.

Using Lemma \ref{lemma: bound on the distance to x^*}, we derive
\begin{align}
    &f(y_{j})- f(x^*) \nonumber\\
    &\leq \delta\left[(1+D/\epsilon) + \left(4C^2 + 2 + D/\epsilon\right) 6 C''/\epsilon + 4\delta\left(C^2 + 1 \right) (1+D/\epsilon) \right].
    % &= \delta \left(1 + \frac{2 L_f/\rho\left[3C''/\epsilon+4\rho^{-1}\left(C + 1 \right) \delta + 2\left(C + 1 \right)/L_f\right]}{\epsilon} \right)\\
    % &\quad + \delta\left(4C + 2 + \frac{2 L_f/\rho\left[3C''/\epsilon+4\rho^{-1}\left(C + 1 \right) \delta + 2\left(C + 1 \right)/L_f\right]}{\epsilon}\right) \frac{6 C''}{\epsilon}\\ 
    % &\quad + 4\delta^2\left(C + 1 \right) \left(1 + \frac{2 L_f/\rho\left[3C''/\epsilon+4\rho^{-1}\left(C + 1 \right) \delta + 2\left(C + 1 \right)/L_f\right]}{\epsilon}\right).
\end{align}
We get $f(y_{j})- f(x^*) = \cO(\delta\epsilon^{-2})$. This concludes the proof of the lemma. 
\end{proof}
\begin{remark}
The conclusion of Lemma \ref{lemma: minimum on xi} is used in its contrapositive form in Theorem \ref{alg: accelerated prox bundle}. Namely, if $y_j$ is not yet $\epsilon$-optimal (i.e., $f(y_j)-f(x^*)>\epsilon$) and iteration $j$ is a null step (i.e., the null step test is not satisfied), then $\xi_j$ is bounded from below as $\xi_j\ge A\epsilon^3$. 
\end{remark}

The following lemma gives an upper bound on $\xi_j$ for serious steps.
\begin{lemma}
\label{lemma: maximum on xi}
There exists a constant $M>0$ such that for all serious step $k$, 
\begin{align}
    \xi_{k+1} \leq \frac{M}{\epsilon^2}. \label{ineq: upper bound on xi}
\end{align}
\end{lemma}
\begin{proof}
We use the smoothness of $f$ to show the claimed bound
\begin{align}
\xi_{k+1} &= f(y_{k+1}) + \frac{\rho}{2}\|y_{k+1}-x_k\|^2 - \left(f_k(y_{k+1}) + \frac{\rho}{2}\|y_{k+1}-x_k\|^2\right) \nonumber\\
&=f(y_{k+1}) - f_k(y_{k+1})\nonumber\\
&\leq f(y_{k+1}) - l_f(y_{k+1}, y_k)\nonumber\\
&\leq \frac{L_f}{2}\|y_{k+1}- y_k\|^2. \label{ineq: bound on first xi in terms of y}
\end{align}
The first equality is by setting $j=k+1$ in the definition \eqref{eq: def of z_j and xi_j}. The first inequality holds because $f_k$ and $f$ have the same value and gradient at $y_k$ and $f_k$ is convex. The last inequality is given by the smoothness of $f$. 

We bound the distance between an iterate and the corresponding proximal center. To do so, we first use the optimality of $y_{k+1}$ for \eqref{eq: bundle minimization alg smooth bundle method} and $f_k\leq f$,
\begin{align*}
    f_k(y_{k+1}) + \frac{\rho}{2}\|y_{k+1} - x_k\|^2 &\leq f_k(x_0) + \frac{\rho}{2}\|x_0 - x_k\|^2 \leq f(x_0) + \frac{\rho}{2}\|x_0 - x_k\|^2.
\end{align*}
Rearranging, 
\begin{align}
    \frac{\rho}{2}\|y_{k+1} - x_k\|^2 &\leq f(x_0) - f_k(y_{k+1}) + \frac{\rho}{2}\|x_0 - x_k\|^2 \nonumber\\
    &\leq f(x_0) - l_f(y_{k+1}, x_0) + \frac{\rho}{2}\|x_0 - x_k\|^2\label{ineq: proof last lemma part 1 1}\\
    &\leq \|\nabla f(x_0)\|\cdot \|x_0 - y_{k+1}\| + \frac{\rho}{2}\|x_0 - x_k\|^2\label{ineq: proof last lemma part 1 2}\\
    &\leq \|\nabla f(x_0)\|\cdot \|y_{k+1} - x_k\| + \|\nabla f(x_0)\| \|x_0 - x_k\| + \frac{\rho}{2}\|x_0 - x_k\|^2\label{ineq: proof last lemma part 1 3}.
\end{align}
Eq. \eqref{ineq: proof last lemma part 1 1} follows from the fact that $f_k \geq l_f(\cdot, x_0)$. Eq. \eqref{ineq: proof last lemma part 1 2} uses the Cauchy-Schwarz inequality, and \eqref{ineq: proof last lemma part 1 3} employs the triangle inequality. Note that \eqref{ineq: proof last lemma part 1 3} is a quadratic inequality in $\|y_{k+1}-x_k\|$. This gives an upper bound
\begin{align}
\|y_{k+1} - x_k\| &\leq \frac{\|\nabla f(x_0)\| + \sqrt{\Delta_k}}{\rho}\notag\\
& \leq \frac{\|\nabla f(x_0)\| + \sqrt{\Delta}}{\rho},\label{ineq: bound on y_{k+1} - x_k}
\end{align}
where $\Delta_k$ and $\Delta$ are defined and related as
\begin{align*}
\Delta_k &:= \|\nabla f(x_0)\|^2 + 2\rho \|\nabla f(x_0)\| \|x_0 - x_k\| + \rho^2 \|x_0 - x_k\|^2\\
& \le \|\nabla f(x_0)\|^2 + 2\rho \|\nabla f(x_0)\|( \|x_0 - x^*\| + \|x^* - x_k\|) + \rho^2 (\|x_0 - x^*\|+\|x^* - x_k\|)^2\\
& \le \Delta := \|\nabla f(x_0)\|^2 + 2\rho \|\nabla f(x_0)\| (\|x_0 - x^*\| + 3C''/\epsilon) + \rho^2 (\|x_0 - x^*\| +  3C''/\epsilon)^2.
\end{align*}  
The  first inequality uses the triangle inquality and the second inequality applies \eqref{ineq: bound on x_k - x^*}.
Now we use this to bound $\|y_{k+1}- y_k\|$ as
\begin{align*}
    \|y_{k+1}- y_k\|&\leq \|y_{k+1}- x_k\| + \|y_{k}- x_{k-1}\| + \|x_{k-1}- x^*\| + \|x^*- x_k\|\\
    &\leq  2\frac{\|\nabla f(x_0)\| + \sqrt{\Delta}}{\rho} + 2\frac{3C''}{\epsilon}.
\end{align*}
% The first inequality uses the triangle inequality, while the second uses \eqref{ineq: bound on y_{k+1} - x_k} and Lemma \ref{lemma: bound on the distance to x^*}. 

Starting from \eqref{ineq: bound on first xi in terms of y},
\begin{align}
    \xi_{k+1} 
&\leq \frac{L_f}{2}\|y_{k+1}- y_k\|^2 \nonumber\\
&\leq \frac{L_f}{2}2\left(2\frac{\|\nabla f(x_0)\| + \sqrt{\Delta}}{\rho}\right)^2 + 4L_f\frac{9C''^2}{\epsilon^2}\label{ineq: final lemma upper bound on xi 1}\\
&\leq 8L_f \frac{\|\nabla f(x_0)\|^2 + \Delta}{\rho^2}+ 4L_f\frac{9C''^2}{\epsilon^2}\label{ineq: final lemma upper bound on xi 2}\\
&\leq \frac{8L_f}{\rho^2} \Bigg( \rho^2 \left(\|x_0 - x^*\| + \frac{3C''}{\epsilon}\right)^2 \label{ineq: final lemma upper bound on xi 3}\\
&\quad + 2\rho \|\nabla f(x_0)\| \left(\frac{3C''}{\epsilon} + \|x_0 - x^*\|\right) + 2\|\nabla f(x_0)\|^2 \Bigg)+ 4L_f\frac{9C''^2}{\epsilon^2}\nonumber\\
&= 8L_f \left(\|x_0 - x^*\| + \frac{3C''}{\epsilon}\right)^2 +  \frac{48 L_f C''\|\nabla f(x_0)\|}{\rho\epsilon} \nonumber\\
&\quad +  \frac{16 L_f}{\rho} \|\nabla f(x_0)\|\|x_0 - x^*\| + \frac{16L_f}{\rho^2}\|\nabla f(x_0)\|^2+ \frac{36 L_f C''^2}{\epsilon^2}\label{ineq: final lemma upper bound on xi 4}\\ 
&= \cO(\frac{1}{\epsilon^2}) \nonumber.
\end{align}

We get \eqref{ineq: final lemma upper bound on xi 1} and \eqref{ineq: final lemma upper bound on xi 2} from the simple inequality \((a+b)^2 \leq 2a^2 + 2b^2\). We apply  \eqref{ineq: bound on x_k - x^*} to get \eqref{ineq: final lemma upper bound on xi 3}. In \eqref{ineq: final lemma upper bound on xi 4}, we simply rearrange the terms.

% This shows that we can define $M$ a constant (independent from $\epsilon$) such that $\xi_k \leq M\epsilon^{-2}$.
\end{proof}
In the above three lemmas, we have shown the following properties of $\xi$ that directly imply a $\cO(\log(\frac{1}{\epsilon}))$ number of consecutive null steps :

\begin{itemize}
    \item $\xi_{j+1} \leq \tau \xi_j$ (with $\tau<1$ independent of $\epsilon$).
    \item $\xi_j \geq A \epsilon^3$ for a null step $j$.
    \item $\xi_{k+1} \leq   M \epsilon^{-2}$ for a serious step $k$.
\end{itemize}
\begin{remark}
    For both the lower bound and the upper bound on $\xi_j$, it is certainly possible to get tighter inequalities (and maybe an upper bound on $\xi$ independent of $\epsilon$). However, this would only affect the constant of the leading term of the convergence rate of at most a factor corresponding to the exponents of the $\epsilon$ in Lemma \ref{lemma: minimum on xi} and Lemma \ref{lemma: maximum on xi}, that is, $3+2 = 5$. 
\end{remark}
\subsection{Overall convergence rate}
\label{section: Overall convergence rate}
The following theorem is the main result of this paper and provides a convergence rate for Algorithm \ref{alg: accelerated prox bundle}. 

\begin{theorem}
\label{thm: iteration complexity}
    Algorithm \ref{alg: accelerated prox bundle} reaches an $\epsilon$-optimal solution in at most the following number of iterations  
    \begin{align*}
       \left(\frac{\sqrt{2\rho}\|x_0-x^*\|}{\sqrt{\epsilon}} +1\right)\left(\frac{5\log(\frac{1}{\epsilon}) + \log\left(\frac{M}{A} \right)}{\log\left(1+ \left( \frac{L_f}{\rho} + \frac{8L_f^2}{\rho^2} + \frac{16L_f^3}{\rho^3} \right)^{-1}\right)}+1\right).
    \end{align*}
\end{theorem}
\begin{proof}
    Let $k$ be the index of a serious step. We have shown in Lemmas \ref{lemma: minimum on xi} and \ref{lemma: maximum on xi} that there exist constants $M$ and $A$, which depend polynomially on the problem parameters and their inverse, such that $\xi_{k+1} \leq M/\epsilon^2$, and for all \textit{null step} index $j > k$ in the sequence of null steps following the serious step $k$, $\xi_{j} > A\epsilon^{3}$. Also, Lemma \ref{lemma: linear convergence of xi} gives the linear convergence of $\xi_{j}$ with rate $\tau = \frac{L_f + \rho C^2}{\rho+L_f + \rho C^2}$. This leads to the following inequality for $j = k+T_k$ in the sequence of null steps,
    \begin{align*}
        &A\epsilon^{3} \leq \xi_j \leq \tau^{T_k-1} \xi_{k+1} \leq \tau^{T_k-1} M/ \epsilon^{2}\\
        \Rightarrow \quad T_k&\leq \frac{5\log(\frac{1}{\epsilon}) + \log\left(\frac{M}{A} \right)}{-\log(\tau)}+1\\
        &= \frac{5\log(\frac{1}{\epsilon}) + \log\left(\frac{M}{A} \right)}{\log\left(1+ \left( \frac{L_f}{\rho} + C^2  \right)^{-1}\right)}+1\\
        &\leq \frac{5\log(\frac{1}{\epsilon}) + \log\left(\frac{M}{A} \right)}{\log\left(1+ \left( \frac{L_f}{\rho} + \frac{8L_f^2}{\rho^2} + \frac{16L_f^3}{\rho^3}  \right)^{-1}\right)}+1.
        % &\leq \frac{5\log(\epsilon^{-1}) - \log(A) +\log(\frac{9}{2}) + \log(L_f)  + 2\log\left(2 + 4\rho\|x_0 - x^*\|^2\right)}{\log\left(1+ \left( \frac{L_f}{\rho} + \frac{16L_f^3}{\rho^3} + \frac{8L_f^2}{\rho^2} + 1 \right)^{-1}\right)}\\
        % &+ \frac{2\log\left(e^2\|x_0 - x^*\|+(\rho \|x_1 - x^*\|^2 + 2(f(x_1) - f(x^*)))\frac{\pi^2}{3}\right)}{\log\left(1+ \left( \frac{L_f}{\rho} + \frac{16L_f^3}{\rho^3} + \frac{8L_f^2}{\rho^2} + 1 \right)^{-1}\right)}+1\\
    \end{align*}
    We now bound the total number of iterations by the maximum number of consecutive null steps times the number of serious steps given in Theorem \ref{thm: accelerated prox point he 1/k2}. This gives the following bound on the total number of iterations
    \begin{align*}
    &\left(\frac{\sqrt{2\rho}\|x_0-x^*\|}{\sqrt{\epsilon}} +1\right)\left(\frac{5\log(\frac{1}{\epsilon}) + \log\left(\frac{M}{A} \right)}{\log\left(1+ \left( \frac{L_f}{\rho} + \frac{8L_f^2}{\rho^2} + \frac{16L_f^3}{\rho^3}  \right)^{-1}\right)}+1\right).
    \end{align*}
\end{proof}

\begin{remark}
    If we choose $\rho = L_f$, we get the following upper bound on the total number of iterations
    \begin{align}
        &\left(\frac{\sqrt{2L_f}\|x_0-x^*\|}{\sqrt{\epsilon}} +1\right)\left(\frac{5\log(\frac{1}{\epsilon}) + \log\left(\frac{M}{A} \right)}{\log\left(1+ \left( 1 + 16 + 8 \right)^{-1}\right)}+1\right)\nonumber\\
        &= \left(\frac{\sqrt{2L_f}\|x_0-x^*\|}{\sqrt{\epsilon}} +1\right)\left(\frac{5}{\log(1+1/25)}\log(\frac{1}{\epsilon}) + K_1\log\left(L_f \right)+K_2\right) \label{rate: total number of iterations rho is Lf},
    \end{align}
    where $K_1$ is a numerical constant and $K_2$ is a constant that depends only on the initial gap in function value and distance to an optimal solution and is independent of $(\epsilon, \rho, L_f)$.    

    In \cite{drori_exact_2017}, it shows that the lower bound for algorithms using only first-order information such as Algorithm \ref{alg: accelerated prox bundle} is asymptotically,
    \begin{align}
        \frac{\sqrt{L_f}\|x_0 - x^*\|}{\sqrt{\epsilon}} \label{lower bound on iteration complexity}.
    \end{align}
    Our rate \eqref{rate: total number of iterations rho is Lf} is suboptimal because it includes an additional $\log(\frac{1}{\epsilon})$ term. The dependence on the smoothness parameter $L_f$ appears under the form $\sqrt{L_f}\log(L_f)$ compared to $\sqrt{L_f}$ in the lower bound \eqref{lower bound on iteration complexity}. 
\end{remark}

% \begin{remark}
%     The previous theorem shows that tuning $\epsilon$ will not result in a change in the complexity which cannot be better than $\cO\left(\epsilon^{-1}\log(\epsilon^{-1})\right)$. \textcolor{blue}{We can certainly get a result of the type "accelerated proximal bundle has an accelerated rate for a large range of parameter value" as stated in \cite{LiangMonteiro2023} and \cite{liang2021proximal}.}
% \end{remark}
\section{Acceleration of the PBM in a nonsmooth setting}
\label{section: Alg for poly objective}
Unlike the gradient descent algorithm, accelerated convergence rates for Nesterov acceleration applied to the proximal point algorithm are not restricted to smooth objective functions. Notably, \cite{Guler1992, He2012, SalzoVilla2012} establish a convergence rate of $\mathcal{O}(\frac{1}{\sqrt{\epsilon}})$ for their algorithms without requiring smoothness of the objective function.

When the proximal subproblem can be solved in a logarithmic number of consecutive null steps, acceleration becomes feasible. This occurs when the objective function is smooth. As demonstrated by \cite{2024proxbundlepolyhedral}, the same holds when dealing with composite problems, combining a smooth term with a polyhedral component \cite{2024proxbundlepolyhedral}. In this section, we prove that applying Nesterov's acceleration to the updates of the proximal center in this framework yields a convergence rate of $\mathcal{O}\left(\frac{1}{\sqrt{\epsilon}}\log\left(\frac{1}{\epsilon}\right)\right)$. This rate, however, depends on the problem's dimension, as does the bound on the number of consecutive null steps.

\subsection{Assumptions}

We study the following composite problem
\begin{align}
    \min_{x \in R^n}\, h(x):=g(x) + f(x).
\end{align}

We assume that $g$ is a convex function and $L_g$-smooth. No longer do we require $f$ to be differentiable; instead, we now assume that $f$ is a convex piecewise linear function. That is, there exists a finite subset $\cV$ of $\R^{n+1}$ such that
\begin{align*}
    \forall x \in \R^{n}, \quad f(x) &= \underset{(v,b) \in \cV}{\max}\, \langle v, x\rangle  + b.
\end{align*}
Let $D$ and $\gamma$ respectively denote the diameter and pyramidal width of $\mathcal{V}$ \cite{lacostejulien2015global, 2024proxbundlepolyhedral}. 

We also define $D_b = \max\{|b_1 - b_2| : (v_1, b_1), (v_2, b_2) \in \cV\}$ and $M_f$ the Lipschitz continuity constant of $f$ such that $M_f = \max\{\|v\| : (v, b) \in \cV\}$.   

We suppose that an oracle can solve the bundle subproblem
\begin{align}
    \min_{y\in R^n}\, g(y)+ f_k(y) + \frac{\rho}{2} \|y - x_j\|^2, \label{eq: composite bundle minimization step}
\end{align}
with $x_j$ the proximal center and $f_k$ a lower model of $f$.

Compared to the previous setting we add the assumption that the set of optimal solutions is bounded.

\subsection{Algorithm description}

The $\epsilon$-subdifferential of $h$ at the point $z \in \R^n$ is the set
\begin{align*}
\partial_{\epsilon} h(z) = \left\{ \xi \in \R^n : h(x) \geq h(z) + \langle \xi, x - z \rangle - \epsilon, \, \forall x \in \R^n \right\},
\end{align*}
and we recall the notation $\operatorname{prox}_{\rho, h}(x) = \argmin_{y\in R^n}\, h(y) + \frac{\rho}{2}\|y - x\|^2$.
\begin{definition}[Definition 2 in \cite{SalzoVilla2012}]
    We say that $y \in \R^n$ is a type 2 approximation of $\operatorname{prox}_{\rho, h}(x_j)$ with $\epsilon_j$-precision and write $y \approx_{2} \operatorname{prox}_{\rho, h}(x_j)$ with $\epsilon_j$-precision if and only if
\begin{align*}
  \rho(x_j - y) \in \partial_{\frac{\epsilon_j^2\rho}{2}} h(y). 
\end{align*}
\end{definition}

We consider the accelerated inexact proximal point algorithm given in \cite{SalzoVilla2012}. It specifies the approximation criterion of Algorithm \ref{alg: proximal point algorithm} using $\approx_2$ and removes the additional sequence $\zeta_j$.
\begin{algorithm}[H]
  \caption{Accelerated Proximal Point Algorithm with Type 2 Approximation}
  \label{alg: accelerated prox point type 2}
\begin{algorithmic}
\Require $h$, $x_0 \in \R^n$
    \State Initialize $\epsilon_0$. Set $y_0 = x_0$,  and $t_0 = 1$ 
    \For{$j \geq 0$}
      \State Compute $y_{j+1} \approx_{2} \operatorname{prox}_{\rho, h}(x_j) \quad$ with $\epsilon_j$-precision
      \State $t_{j+1} \gets \frac{1+\sqrt{1 + 4t_j^2}}{2}$
      \State $x_{j+1} \gets y_{j+1} + \frac{t_j - 1}{t_{j+1}}(y_{j+1} - y_j)$
      \State Update $\epsilon_{j+1}$
    \EndFor
\end{algorithmic}
\end{algorithm}
Using a piecewise linear model such as in the classic PBM (Algorithm \ref{alg: prox bundle}), we derive the following accelerated PBM.
\begin{algorithm}[htb]
  \caption{Accelerated Proximal Bundle Method}
  \label{alg: accelerated prox bundle for polyhedral}
\begin{algorithmic}
\Require $f$, $g$, $x_0 \in \R^n$, $\mathcal{I}_0$ an index set of initial cuts, $B>0$ and $\rho>0$
    \State Set $y_0 = \zeta_0 = x_0$, $t_0 = 1$ and $j = 0$
    \For{$k \geq 0$}
      \State Compute $y_{k+1}$ solving the following convex QCQP
      \begin{align}
             & \underset{t \in \mathbb{R}, y \in \mathbb{R}^n}{\min} \; t + g(y) + \frac{\rho}{2} \|y - x_j\|^2 \label{eq: bundle minimization alg polyhedral bundle method}\\
    \text{s.t.} \quad &\forall i \in \mathcal{I}_k, \quad t - f(y_i) - \langle v_i, y - y_i\rangle 
    \geq 0. \nonumber 
      \end{align}
        \State Compute $f(y_{k+1})$ and $v_{k+1} \in \partial f(y_{k+1})$
      \State $\mathcal{I}_{k+1} \leftarrow \mathcal{I}_k \cup \{k+1\}$
      % \State \textbf{//Null step test} : 
      \If{$f(y_{k+1}) - f_{k}(y_{k+1}) \leq \frac{\rho \epsilon_j^2}{2}$} \Comment{(null step test)}
        \State $t_{j+1} \gets \frac{1+\sqrt{1 + 4t_j^2}}{2}$
        \State $\zeta_{j+1} \gets y_{k+1}$
        \State $x_{j+1} \gets \zeta_{j+1} + \frac{t_j - 1}{t_{j+1}}(\zeta_{j+1} - \zeta_{j})$ \Comment{(serious step)}
        \State  $\epsilon_{j+1} \gets \frac{\sqrt{6B}}{\pi\sqrt{\rho}(j+3)^2}$
        \State $j \gets j+1$
      \EndIf
      % \State $(\hat{v}_{k+1}, \hat{b}_{k+1}) \in \underset{(v, b) \in \mathcal{V}}{\argmax}\; (v^{T} y_{k+1} + b)$
    \EndFor
\end{algorithmic}
\end{algorithm}

\paragraph{Comparison between Algorithm \ref{alg: accelerated prox bundle} and Algorithm \ref{alg: accelerated prox bundle for polyhedral}} Both algorithms share a similar structure and can be viewed as accelerated versions of the classic PBM (Algorithm \ref{alg: prox bundle}). The Nesterov acceleration is implemented using the same sequence, $t_k$. Additionally, the bundle minimization step \eqref{eq: bundle minimization alg polyhedral bundle method} in Algorithm \ref{alg: accelerated prox bundle for polyhedral} can be interpreted as the bundle minimization step \eqref{eq: bundle minimization alg smooth bundle method} in Algorithm \ref{alg: accelerated prox bundle}, where the smoothness constant $L_f = +\infty$ since $f$ is not assumed to be differentiable.

Despite these similarities, the two algorithms differ in their null step tests. In Algorithm \ref{alg: accelerated prox bundle for polyhedral}, the test requires knowledge of the desired accuracy, $\epsilon$. Furthermore, Algorithm \ref{alg: accelerated prox bundle} keeps track of the sequence $\zeta$ which is updated as $\zeta_{k+1} \gets x_k - \frac{1}{\rho}\nabla f(y_{k+1})$. 

\subsection{Analysis of the algorithm}
\subsubsection{Number of serious steps}
The following proposition establishes a connection between the two criteria: the gap between the value of the lower model and the true function, as typically considered in the cutting-plane literature, and the inexactness criterion $\approx_2$ from the proximal point algorithm literature.
\begin{proposition}
\label{prop: type 2 approx is satisfied}
    Let $y_{k+1}$ be an iterate of the Algorithm \ref{alg: accelerated prox bundle for polyhedral} such that the null step test $f(y_{k+1}) - f_k(y_{k+1}) \leq \frac{\rho \epsilon_j^2}{2}$ is satisfied for some $j\geq 0$. Then, 
    \begin{align*}
        y_{k+1} \approx_2 \operatorname{prox}_{\rho, h}(x_j) \quad \text{with $\epsilon_j$-precision}.
    \end{align*}
\end{proposition}
\begin{proof}
    Let $u \in R^n$. 
    \begin{align*}
        g(u) + f(u) + \frac{\rho}{2} \|u - x_j\|^2 &\geq  g(u) + f_k(u) + \frac{\rho}{2} \|u - x_j\|^2\\
        &\geq g(y_{k+1}) + f_k(y_{k+1}) + \frac{\rho}{2}\|y_{k+1} - x_j\|^2 + \frac{\rho}{2}\|y_{k+1} - u\|^2.\\
        &\geq g(y_{k+1}) + f(y_{k+1}) + \frac{\rho}{2}\|y_{k+1} - x_j\|^2 + \frac{\rho}{2}\|y_{k+1} - u\|^2 - \frac{\rho \epsilon_j^2}{2}.
    \end{align*}
    The first line inequality uses $f \geq f_k$. The second inequality leverages the optimality of $y_{k+1}$ for the bundle minimization step \eqref{eq: composite bundle minimization step} and the strong convexity of this problem. The last inequality uses the hypothesis that $f(y_{k+1}) - f_k(y_{k+1}) \leq \frac{\rho \epsilon_j^2}{2}$.

    Rearranging, we get 
    \begin{align*}
        g(u) + f(u)  &\geq g(y_{k+1}) + f(y_{k+1}) + \frac{\rho}{2} \left(\|y_{k+1} - x_j\|^2 + \|y_{k+1} - u\|^2 -  \|u - x_j\|^2 \right) - \frac{\rho \epsilon_j^2}{2} \\
        &= g(y_{k+1}) + f(y_{k+1}) + \rho \langle u -y_{k+1}, x_j -y_{k+1}\rangle - \frac{\rho \epsilon_j^2}{2}.
    \end{align*}
    As this inequality holds for all $u \in \R^n$, this proves that $\rho(x_j - y_{k+1}) \in \partial_{\frac{\epsilon_j^2\rho}{2}} h(y_{k+1})$. Following the definition of $\approx_{2}$ with $\epsilon_j$-precision, this concludes the proof.
\end{proof}

% The previous Proposition \ref{prop: type 2 approx is satisfied} has the following a reciprocal. 
% \begin{proposition}
%     Let $x_j$ a proximal center and suppose that $y_{k+1}$ is such that 
%     \begin{align*}
%         y_{k+1} \approx_2 \operatorname{prox}_{\rho, h}(x_jh) \quad \text{with $\epsilon_j$-precision}.
%     \end{align*}
%     Then, 
%     \begin{align*}
%         f(y_{k+1}) - f_k(y_{k+1}) \leq SOMETHING.
%     \end{align*}
% \end{proposition}
% \begin{proof}
%     Because $y_{k+1}$ approximates the proximal problem with $\epsilon_j$-precision, noting $\eta = \frac{\rho \epsilon_j^2}{2}$
%     \begin{align}
%         \rho(x_j - y_{k+1}) &\in \partial_\eta (g+f) (y_{k+1}) \nonumber\\
%         &= \bigcup_{\eta_1 + \eta_2 = \eta} \partial_{\eta_1} g (y_{k+1}) + \partial_{\eta_2} f (y_{k+1}) \\
%         &\subseteq  \partial_\eta g (y_{k+1}) + \partial_\eta f (y_{k+1})\\
%         &\subseteq \nabla g (y_{k+1}) + B(0, \sqrt{2L_g\eta}) + \partial_\eta f (y_{k+1}) \label{inclusion: reciprocal of type 2 approx}\\
%     \end{align}
% Following Lemma C.2 in \cite{2024proxbundlepolyhedral} we define $w_{k+1} = \nabla (g (y_{k+1}) + \frac{\rho}{2}\|y_{k+1}-x_j\|^2$. The inclusion \eqref{inclusion: reciprocal of type 2 approx} translates into 
% \begin{align}
%     w_{k+1} \in \left(\partial_\eta f (y_{k+1}) \right) + B(0, \sqrt{2L_g\eta}).
% \end{align}

% \end{proof}
\begin{corollary}
    Algorithm \ref{alg: accelerated prox bundle for polyhedral} is an instance of Algorithm \ref{alg: accelerated prox point type 2} with $h = f + g$. 
\end{corollary}
\begin{proof}
    After noticing that Algorithm \ref{alg: accelerated prox bundle for polyhedral} contains an implicit loop in which an approximate solution of the proximal problem is computed and during which the proximal center is not updated, this corollary follows directly from Proposition \ref{prop: type 2 approx is satisfied}.
\end{proof}
\begin{theorem}[Number of serious steps of Algorithm \ref{alg: accelerated prox bundle for polyhedral}]
\label{Number of serious steps of PPM type 2}
    At the $j^{th}$ serious step, the iterate $\zeta_j$ is such that
    \begin{align}
        h(\zeta_j) - h^* \leq \beta_j (h(x_0) - h^* + \frac{\rho}{2} \|x_0 - x^*\|^2) + \delta_j \label{ineq: bound on gap for type 2 approx PPA}
    \end{align}
with $\beta_j \leq \frac{1}{\left(1 + \frac{j}{2} \right)^2}$ and $\delta_j \leq \frac{\beta_j}{2} \sum_{i=0}^{j-1} \rho \epsilon_i^2 \left(1 + \sqrt{2}(i+1) \right)^2.$

Moreover, for $B>0$ choosing $\epsilon_j \leq \frac{\sqrt{2B}\cdot \sqrt{6}/\pi}{\sqrt{\rho}(j+1) \left(1 + \sqrt{2}(j+1)\right)}$ leads to $\delta_j \leq B \beta_j$. 
We can then rewrite the inequality \eqref{ineq: bound on gap for type 2 approx PPA} as
    \begin{align}
        h(\zeta_j) - h^* \leq \frac{4}{\left(j+2\right)^2} (h(x_0) - h^* + \frac{\rho}{2} \|x_0 - x^*\|^2 + B). \label{ineq: nb serious steps in terms of j}
    \end{align}
\end{theorem}
\begin{proof}
    A more general proof is given in Remark 4.7 of \cite{SalzoVilla2012}. Compared to their notations, we have $\forall i, \lambda_i = \rho^{-1}$, $a = 1$, $A_0 = \rho$, and $p=1$. 
\end{proof}
This theorem shows that the number of serious steps to reach an $\epsilon$-optimal solution is at most 
\begin{align}
    \frac{2\sqrt{h(x_0) - h^* + \frac{\rho}{2} \|x_0 - x^*\|^2 + B}}{\sqrt{\epsilon}}. \label{number of serious step polyhedral} 
\end{align}
    
\begin{lemma}
\label{lemma: lower bound on epsilonj}
    While an $\epsilon$-optimal solution has not been found, choosing $\epsilon_j = \frac{\sqrt{6B}}{\pi\sqrt{\rho}(j+2)^2}$ satisfies the upper bound on $\epsilon_j$ in Theorem \ref{Number of serious steps of PPM type 2}; moreover, such $\epsilon_j$ is lower bounded as 
    \begin{align}
        \forall j, \quad  \epsilon_j \geq \frac{\sqrt{6B}\epsilon}{4\pi\sqrt{\rho}\left(h(x_0) - h^* + \frac{\rho}{2} \|x_0 - x^*\|^2 + B\right)}.
    \end{align}
\end{lemma}

\begin{proof}
    We suppose that $h(\zeta_j) - h^* \geq \epsilon$. Following \eqref{ineq: nb serious steps in terms of j}, 
    \begin{align*}
       \epsilon \leq \frac{4}{\left(j+2\right)^2} (h(x_0) - h^* + \frac{\rho}{2} \|x_0 - x^*\|^2 + B).
    \end{align*}
    Rearranging and plugging the maximum value of $j$ (given in \eqref{number of serious step polyhedral}) in the expression of $\epsilon_j$ leads to the claimed bound.
\end{proof}
\subsubsection{Number of consecutive null steps}

\begin{theorem}[Lemma 3.6 in \cite{2024proxbundlepolyhedral}]
\label{theorem: Lemma 3 in our paper}
In Algorithm \ref{alg: accelerated prox bundle for polyhedral}, at most the following number of consecutive null steps occur between two serious steps 
\begin{align}
    \left[1+\max\left\{2,\frac{D^2}{{\mu}_{\psi, \rho, j} \rho \gamma^2}\right\} \log \left(\frac{4D^4}{\rho^4\epsilon_j^2} \right)\right],
\end{align}
with ${\mu}_{\psi, \rho}$ such that 
\begin{align}
    \frac{1}{2}{\mu}_{\psi, \rho, j}^{-1} &= D_b + 3 \frac{8 M_f^2}{\rho}  + 6M_f \|x^*_j\| + 2L_g \left[\left(\frac{4M_f}{\rho} + \|x^*_j\|\right)^2 + 1\right], \label{def: mu psi rho j} 
\end{align}
where $x_j^*=\argmin_x h(x)+\frac{\rho}{2}\|x-x_j\|^2$.
\end{theorem}

\begin{definition}[Level-Boundedness]
A function $f: \mathbb{R}^n \to \mathbb{R}$ is said to be \emph{level-bounded} if for every $c \in \mathbb{R}$, the level set
\begin{align*}
    L_c = \{ x \in \mathbb{R}^n : f(x) \leq c \}
\end{align*}
is bounded.
\end{definition}

We recall the following known result (Lemma \ref{lm: level boundedness}) that derives the level-boundedness of the convex function $h$ based on the boundedness of its minimizers. For completeness, a proof is given in the Appendix \ref{sec: proof of level-boundedness}. 
\begin{lemma}[Level-Boundedness of Convex Functions with Bounded Minimizers]\label{lm: level boundedness}
    If $\{x \in \R^n : f(x) =   \min_{y \in \R^n}\, f(y)\}$ is non-empty and bounded, then $f$ is level bounded. 
\end{lemma}
% \begin{proof}
% Since the set of minimizers of $f$ is non-empty, the minimum $c_1 := \min_{y \in \R^n}\, f(y)$ is finite.  
% We proceed by contradiction. By the assumption of the lemma, the level set $L_1 = \{ x : f(x) \leq c_1 \}$ is bounded. Suppose there exists $c_2>c_1$ such that $L_2 = \{ x : f(x) \leq c_2 \}$ is unbounded. Note that $L_1 \subseteq L_2$. Since $L_2$ is unbounded, it follows that the recession cone of $L_2$, denoted by $0^+ L_2$, is nontrivial; that is, there exists a nonzero direction $d \in 0^+ L_2$ (see, for instance, \cite{rockafellar1970convex} Theorem 8.4).

% In particular, for any $x_1 \in L_1$, we have $x_1 + td \in L_2$ for all $t \geq 0$. Define the function $\varphi(t) = f(x_1 + td)$. Since $x_1 + td \in L_2$ for all $t \geq 0$, $\varphi(t)$ must be nonincreasing. Otherwise, for some $t > 0$, we would have $f(x_1 + td) > c_2$, which contradicts the assumption that $x_1 + td \in L_2$ for all $t\ge 0$.

% Moreover, since $\varphi(t)$ is nonincreasing and $f(x_1) \leq c_1$, it follows that $f(x_1 + td) \leq c_1$ for all $t \geq 0$, implying that $x_1 + td \in L_1$ for all $t \geq 0$. This contradicts the boundedness of $L_1$. Therefore, all level sets of $f$ must be bounded.
% \end{proof}

The constant ${\mu}_{\psi, \rho, j}$ defined in Theorem \ref{theorem: Lemma 3 in our paper} depends on $\|x^*_j\|$, the norm of the proximal solution $\operatorname{prox}_{\rho, h}(x_j)$. We leverage the level-boundedness of $h$ to bound $\|x^*_j\|$ uniformly in $j$. 
\begin{proposition}
    The optimal solutions of the proximal problems $(x_j^*)_{j\in \N}$ are bounded as
    \begin{align}
        \forall j,\quad \|x_j^*\| \leq 3\|x^*\| + 6 R_{h, x_0},
    \end{align}
    where $R_{h, x_0}$ is the radius of a Euclidean ball centered at zero that contains the level set of $h$ at level $h(x_0) + \frac{\rho}{2} \|x_0 - x^*\|^2 + B$. Note that $R_{h,x_0}$ does not depend on $\epsilon$.
    % depends on $h$ and $x_0$ but not on $\epsilon$. 
\end{proposition}
\begin{proof}
By Theorem \ref{Number of serious steps of PPM type 2}, $\beta_j\le 1$ for all $j$. Following \eqref{ineq: nb serious steps in terms of j}, all the $\zeta_j$ are in the level set of $h$ with level $h(x_0) + \frac{\rho}{2} \|x_0 - x^*\|^2 + B$. 

    We combine this bound with the update of the proximal center given in Algorithm \ref{alg: accelerated prox bundle for polyhedral}
    \begin{align*}
        \|x_{j+1}\| &= \left\|\zeta_{j+1} + \frac{t_j - 1}{t_{j+1}}(\zeta_{j+1} - \zeta_{j}) \right\|\\
        &\leq \|\zeta_{j+1}\| + \frac{t_j - 1}{t_{j+1}}(\|\zeta_{j+1}\| + \|\zeta_{j}\|)\\
        &\leq 3 R_{h, x_0}.
    \end{align*}
    For the last inequality we use $\frac{t_j - 1}{t_{j+1}} \leq 1$ as in \eqref{ineq: ratio of t smaller than 1}. 
    
    By optimality of $x^*$ and $x_j^*$ for their respective problems,
\begin{align*}
    h(x^*) + \frac{\rho}{2}\|x^*_j -x_{j}\|^2 \leq h(x^*_j) + \frac{\rho}{2}\|x^*_j -x_{j}\|^2 \leq h(x^*) + \frac{\rho}{2}\|x^* -x_{j}\|^2.
\end{align*}

Thus, $\|x^*_j -x_{j}\| \leq \|x^* -x_{j}\|$. We combine these two inequalities,
\begin{align}
    \|x_j^*\| &\leq \|x^*\| + \|x_j^* - x^*\| \leq  \|x^*\| + \|x_j^* - x_{j}\| + \|x^*- x_{j}\| \nonumber\\
    &\leq  \|x^*\| + 2\|x^*- x_j\| \leq 3\|x^*\| + 2\|x_j\| \nonumber\\
    &\leq  3\|x^*\| + 6 R_{h, x_0} \label{ineq: uniform bound on xj star}.
\end{align}
\end{proof}

\begin{theorem}
    Algorithm \ref{alg: accelerated prox bundle for polyhedral} provides an $\epsilon$-optimal solution in at most the following number of iterations 
    \begin{align}
        \frac{2\sqrt{h(x_0) - h^* + \frac{\rho}{2} \|x_0 - x^*\|^2 + B}}{\sqrt{\epsilon}} \cdot \left[1+\max\left\{2,\frac{D^2}{{\mu}_{\psi, \rho} \rho \gamma^2}\right\} \log \left(\frac{4D^4}{\rho^4 \sigma \epsilon^2} \right)\right]
    \end{align}
    with
    \begin{align*}
        \frac{1}{2}\mu_{\psi, \rho}^{-1}  \defeq D_b {+} 3 \frac{8 M_f^2}{\rho}  {+} 6M_f (3\|x^*\| {+} 6 R_{h, x_0}) {+} 2L_g \left[\left(\frac{4M_f}{\rho} {+} 3\|x^*\| {+} 6 R_{h, x_0}\right)^2 {+} 1\right],
    \end{align*}
    and 
    \begin{align}
        \sigma = \frac{{6}{B}}{\pi^2{\rho}\left(h(x_0) - h^* + \frac{\rho}{2} \|x_0 - x^*\|^2 + B\right)^2}.
    \end{align}
\end{theorem}

\begin{proof}
We upper bound the total number of iterations by the bound on the number of serious steps given in \eqref{number of serious step polyhedral} and the bound on the number of consecutive null steps given in Theorem \ref{theorem: Lemma 3 in our paper}. 

Combining the definition of $\mu_{\psi, \rho, j}$ in \eqref{def: mu psi rho j} and \eqref{ineq: uniform bound on xj star}, we get $\forall j, \, \mu_{\psi, \rho} \leq \mu_{\psi, \rho, j}$.

Lemma \ref{lemma: lower bound on epsilonj} shows that for all serious step index $j$ such that an $\epsilon$-optimal solution has not been found yet, $\sigma \epsilon^2 \leq \epsilon_j^2$.

\end{proof}

\section{Remarks and conclusion}
\paragraph{Bundle management and varying $\rho$ parameter}
For simplicity, throughout this article, we have assumed that the parameter $\rho$ remains fixed and that all cuts are retained in memory during sequences of null steps. However, the proposed algorithms can be readily adapted to more general settings where $\rho$ may be updated at serious steps following any nonincreasing sequence as Theorem \ref{thm: accelerated prox point he 1/k2} still holds. Also, selective cut management strategies may be employed, without significant changes to the convergence properties. 

\paragraph{Conclusion}
In this paper, we have proposed an Accelerated Proximal Bundle algorithm designed to improve the convergence rate of the proximal bundle methods for smooth objectives. By incorporating a novel null step test and a smooth lower model, we demonstrated that our algorithm achieves an improved convergence rate of \(\cO(\frac{1}{\sqrt{\epsilon}} \log(\frac{1}{\epsilon}))\). This enhancement addresses an open question about achieving accelerated rates with proximal bundle methods. We have further proposed an acceleration of the proximal bundle method for composite convex optimization with piecewise linear nonsmoothness. The resulting Algorithm 5 achieves a similar rate of \(\cO(\frac{1}{\sqrt{\epsilon}} \log(\frac{1}{\epsilon}))\), although the complexity bound also depends the problem's dimension through some geometric quantity such as the pyramidal width of the subdifferentials of the piecewise linear function. 
% This complexity bound improves on the authors' recent work \cite{2024proxbundlepolyhedral}.

% \section*{Ethics declaration}
% \subsection*{Conflict of interest}
% The authors have no conflict of interest to declare that are relevant to the content of this article.

\backmatter

% \bmhead{Supplementary information}

% If your article has accompanying supplementary file/s please state so here. 

% Authors reporting data from electrophoretic gels and blots should supply the full unprocessed scans for key as part of their Supplementary information. This may be requested by the editorial team/s if it is missing.

% Please refer to Journal-level guidance for any specific requirements.

% \bmhead{Acknowledgements}

\begin{appendices}
\section{Additional Proofs}
\subsection{Proof of Proposition \ref{proposition: can find a sequence converging to the sup}} \label{sec: proof of prop 12}
\begin{proof}
    Let $x_\delta \in \R^n$ such that $\bar{h} - h(x_\delta) \leq \delta$. 
    We consider the problem
    \begin{align}
        \max_{x \in \R^n}\, h(x) - \frac{\delta}{2}\|x - x_\delta\| \label{pb: proof for vanishing gradient sequence}.
    \end{align}
    As $h$ is upper-bounded, \eqref{pb: proof for vanishing gradient sequence} can be restricted to a compact feasible set. By continuity of the objective, this problem admits an optimal solution $x_\delta^*$. In particular, as $h(x_\delta^*) \geq h(x_\delta)$, $\|\nabla h(x_\delta^*)\|\leq \delta$ would conclude the proof. We assume the opposite.  $h$ is differentiable at $x_\delta^*$ so there exists $t>0$ such that $z = x_\delta^* + t\nabla h(x_\delta^*)$ satisfies:
    \begin{align*}
        h(z) = h(x_\delta^* + t\nabla h(x_\delta^*)) &\geq h(x_\delta^*) + \frac{3}{4} \left(t\nabla h(x_\delta^*)^T\right) \nabla h(x_\delta^*) \\
        &= h(x_\delta^*) + \frac{3}{4}\|x_\delta^* - z\|\cdot \|\nabla h(x_\delta^*)\| \\
       \implies\, h(z) - \frac{\delta}{2}\|z - x_\delta\| &\geq h(x_\delta^*) + \delta \frac{3}{4}\|x_\delta^* - z\| - \frac{\delta}{2}\|z - x_\delta\| \\
       &\geq h(x_\delta^*) + \delta \frac{3}{4}\|x_\delta^* - z\| - \frac{\delta}{2}(\|x_\delta^* - z\| + \|x_\delta^* - x_\delta\|) \\
       &\geq h(x_\delta^*)  - \frac{\delta}{2}\|x_\delta^* - x_\delta\| + \frac{\delta}{4}\|x_\delta^* - z\|.
    \end{align*}
    This contradicts the maximality of $x_\delta^*$ and concludes the proof.
\end{proof}

\subsection{Proof of Lemma \ref{lm: recurrence inequality}} \label{sec: proof of Lemma 18}

\begin{proof}
We first consider the homogeneous recurrence $s_0 = r_0$ and
\begin{align*}
    s_{k+1} &= (1 + \frac{2}{k+2})s_{k}\\
    \Leftrightarrow \log(s_{k+1}) - \log(s_{k}) &= \log(1 + \frac{2}{k+2})\\
    &\leq \frac{2}{k+2}. &&\text{(concavity of the log)}
\end{align*}
We sum up these inequalities,
\begin{align*}
    \log(s_{k+1}) - \log(s_{0}) &\leq 2\sum_{i=0}^k \frac{1}{i+2}\\
    &\leq 2\sum_{i=1}^{k+1} \frac{1}{i}\\
    &\leq 2+ 2\log(k+1).
\end{align*}
This yields : $s_{k+1} \leq s_0\cdot e^2\cdot (k+1)^2$, for $k\geq 0$. 

Using the inequality $\forall x>-1 , \,\log(1+x) \geq x - x^2/2$ : 
\begin{align*}
    \log(s_{k+1}) - \log(s_{0}) &\geq 2\sum_{i=0}^k \frac{1}{i+2} - \frac{1}{2}\sum_{i=0}^k \left(\frac{2}{i+2} \right)^2\\
    &=2\sum_{i=2}^{k+2} \frac{1}{i} - 2\sum_{i=2}^{k+2} \left(\frac{1}{i} \right)^2\\
    &\geq -1 + 2\log(k+3) -2\frac{\pi^2}{6}.
\end{align*}
This gives $s_{k} \geq s_0 \exp({-1-\frac{\pi^2}{3}})(k+2)^2$ for $k>0$.
We now show that the inhomogeneous term $\frac{2C'}{k+2}$ vanishes at infinity.

For $k\ge 0$, let $w_{k} = \frac{r_{k}}{s_{k}}$. We have $w_0 = 1$. 
\begin{align*}
    w_{k+1} &=  \frac{r_{k+1}}{s_{k+1}} 
    \leq \frac{(1+\frac{2}{k+2}) r_k + \frac{2C'}{k+2}}{(1 + \frac{2}{k+2})s_{k}}
    = w_{k} + \frac{\frac{2C'}{k+2}}{(1 + \frac{2}{k+2})s_{k}}
    = w_{k} + \frac{2C'}{(k+4)s_{k}}\\
    &\leq w_{k} + \frac{2C'}{(k+4) s_0 \exp({-1-\frac{\pi^2}{3}})(k+2)^2} 
    \leq w_{k} + \frac{2C'}{s_0\cdot \exp({-1-\frac{\pi^2}{3}})\cdot (k+1)^3}.
\end{align*}
Solving this simple linear recurrence leads to $\forall k\ge 0$,
\begin{align*}
    w_{k}&\leq w_0+\frac{2C'}{s_0\cdot \exp({-1-\frac{\pi^2}{3}})}\cdot \sum_{i=1}^{+\infty}{\frac{1}{k^3}}
    \leq 1+\frac{2C'\pi^2}{6 s_0\cdot \exp({-1-\frac{\pi^2}{3}})}.
\end{align*}

Combining the bounds on $w_k$ and $s_k$, we get
\begin{align*}
    \forall k\ge 1, \quad r_k &= w_{k} s_k
    \leq \left(e^2\cdot r_0+\frac{C'\pi^2e^{({3+\frac{\pi^2}{3}})}}{3}\right) \cdot k^2.
\end{align*}
\end{proof}

\subsection{Proof of Lemma \ref{lm: level boundedness}}\label{sec: proof of level-boundedness}
\begin{proof}
Since the set of minimizers of $f$ is non-empty, the minimum $c_1 := \min_{y \in \R^n}\, f(y)$ is finite.  
We proceed by contradiction. By the assumption of the lemma, the level set $L_1 = \{ x : f(x) \leq c_1 \}$ is bounded. Suppose there exists $c_2>c_1$ such that $L_2 = \{ x : f(x) \leq c_2 \}$ is unbounded. Note that $L_1 \subseteq L_2$. Since $L_2$ is unbounded, it follows that the recession cone of $L_2$, denoted by $0^+ L_2$, is nontrivial; that is, there exists a nonzero direction $d \in 0^+ L_2$ (see, for instance, \cite{rockafellar1970convex} Theorem 8.4).

In particular, for any $x_1 \in L_1$, we have $x_1 + td \in L_2$ for all $t \geq 0$. Define the function $\varphi(t) = f(x_1 + td)$. Since $x_1 + td \in L_2$ for all $t \geq 0$, $\varphi(t)$ must be nonincreasing. Otherwise, for some $t > 0$, we would have $f(x_1 + td) > c_2$, which contradicts the assumption that $x_1 + td \in L_2$ for all $t\ge 0$.

Moreover, since $\varphi(t)$ is nonincreasing and $f(x_1) \leq c_1$, it follows that $f(x_1 + td) \leq c_1$ for all $t \geq 0$, implying that $x_1 + td \in L_1$ for all $t \geq 0$. This contradicts the boundedness of $L_1$. Therefore, all level sets of $f$ must be bounded.
\end{proof}
% \section{Section title of first appendix}\label{secA1}

% \newpage




\end{appendices}


\bibliography{sn-bibliography}% common bib file

\end{document}